\documentclass[11pt]{amsart}
\usepackage[margin=1.03in]{geometry}
\usepackage{amsmath,amssymb,mathtools}
\usepackage{microtype}
\usepackage[hidelinks]{hyperref}
\usepackage{booktabs}
\newtheorem{theorem}{Theorem}[section]
\newtheorem{proposition}[theorem]{Proposition}
\newtheorem{lemma}[theorem]{Lemma}
\newtheorem{corollary}[theorem]{Corollary}
\newtheorem{remark}[theorem]{Remark}
\newcommand{\Ric}{\operatorname{Ric}}
\newcommand{\ind}{\operatorname{ind}}

\newcommand{\F}{\mathbb F}
\title[Morse index under Hopf squashing]{Morse index and bifurcation of Clifford hypersurfaces under Hopf squashing}
\author{Anna Siffert}
\date{}

\newcommand{\Oct}{\mathbb O}
\newcommand{\Sph}{\mathbb S}
\newcommand{\HH}{\mathcal H}
\newcommand{\PP}{\mathcal P}
\newcommand{\VV}{\mathcal V}
\newcommand{\av}[1]{\langle #1\rangle}
\begin{document}

\begin{abstract}
Clifford hypersurfaces are products of two spheres embedded in a round sphere. We call them balanced when their radii are chosen so that the product is minimal in the round metric. We study how their Morse index and nullity change when the ambient metric is deformed along the fibers of a Hopf fibration. This deformation, called \emph{Hopf squashing}, rescales the fiber metric by a parameter $t>0$ while leaving the horizontal metric unchanged. For the complex, quaternionic, and octonionic Hopf fibrations, we determine the index and nullity of balanced Clifford hypersurfaces throughout the fiber-shortening range and show that their Jacobi potential is independent of $t$. In the octonionic case we also determine the first positive Jacobi eigenvalue and its multiplicity. In the octonionic case, at $t=27$, a branch of nonhomogeneous minimal hypersurfaces bifurcates from the balanced Clifford hypersurface; its members have Morse index $1432$ and nullity $15$.
\end{abstract}

\maketitle

\section{Introduction}
We consider the complex and quaternionic Hopf fibrations
\[
 S^1\longrightarrow S^{2N-1}\longrightarrow\mathbb CP^{N-1},\qquad
 S^3\longrightarrow S^{4N-1}\longrightarrow\mathbb HP^{N-1},\qquad N\ge2,
\]
and the octonionic fibration
\[
 S^7\longrightarrow S^{15}\xrightarrow{\,\pi\,}S^8(1/2).
\]
In each case, let $\mathcal V$ be the tangent distribution to the Hopf fibers and let $\mathcal H=\mathcal V^\perp$ be its orthogonal complement for the round metric $g_1$. Writing $g_{\mathcal H}$ and $g_{\mathcal V}$ for the corresponding parts of the round metric, we set
\[
 g_t=g_{\mathcal H}+t\,g_{\mathcal V},\qquad t>0.
\]
We call this deformation \emph{Hopf squashing}. The parameter $t$ is the squared fiber-length scale: $t=1$ gives the round metric, $0<t<1$ shortens the fibers, and $t>1$ stretches them.

%kontrolle

For a uniform statement covering all three fibrations, let $\mathbb F=\mathbb C,\mathbb H,$ or $\mathbb O$, set $d=\dim_{\mathbb R}\mathbb F$ and $r=d-1$. In the complex and quaternionic cases write $N=p+q+2$, with $p,q\ge0$; in the octonionic case take $p=q=0$ and $N=2$.
Define
\[
 a=d(p+1)-1,\qquad b=d(q+1)-1,\qquad m=a+b=dN-2
\]
and the balanced Clifford hypersurface
\begin{equation}\label{eq:Sigma}
 \Sigma_{p,q}=S^a\!\left(\sqrt{\frac am}\right)
 \times S^b\!\left(\sqrt{\frac bm}\right)\subset S^{m+1}.
\end{equation}
Let $\nu$ denote a unit normal field along $\Sigma_{p,q}$ with respect to $g_t$, and let $A_t$ be its shape operator, defined by $A_t(X)=-\nabla^{g_t}_X\nu$ for tangent vectors $X$ to $\Sigma_{p,q}$. We write $|A_t|^2$ for the sum of the squares of its principal curvatures, and $\operatorname{Ric}_{g_t}(\nu,\nu)$ for the Ricci curvature in the normal direction.

We use the nonnegative Laplace operator $\Delta_{\Sigma,g_t}=-\operatorname{div}_{\Sigma,g_t}\nabla_{\Sigma,g_t}$. The Jacobi operator governing the second variation of area (see, for example, \cite{Simons}) is
\[
 J_t=\Delta_{\Sigma,g_t}-\bigl(|A_t|^2+\operatorname{Ric}_{g_t}(\nu,\nu)\bigr).
\]
The function $|A_t|^2+\operatorname{Ric}_{g_t}(\nu,\nu)$ is called the \emph{Jacobi potential}. Throughout, the Laplacian appearing in the Jacobi operator is nonnegative. In the bifurcation analysis below, the symbols $\Delta_X$ and $\Delta_Y$ instead denote the geometric Laplacians $\operatorname{div}\nabla$ on the unit spheres, with nonpositive eigenvalues; likewise $L_{\VV}=-\Delta_{\VV}$ is nonnegative.
The spectrum of $J_t$ consists of its eigenvalues, counted with multiplicity. The Morse index is the total multiplicity of its negative eigenvalues, and the nullity is the dimension of its zero eigenspace. 

For the balanced Clifford hypersurfaces, the squared norm of the second fundamental form and the ambient normal Ricci curvature satisfy
\[
 |A_t|^2+\operatorname{Ric}_{g_t}(\nu,\nu)=2m.
\]
Consequently, the Jacobi operator changes with $t$ only through the induced Laplacian. We determine the index and nullity in the fiber-shortening range and the first positive Jacobi eigenvalue in the octonionic case. We also establish bifurcation at $t=27$ in the octonionic family and compute the index and nullity along the resulting nonhomogeneous branch.

The following theorem describes the index and nullity for all three Hopf fibrations. Here $\operatorname{ind}$ and $\operatorname{nul}$ refer to index and nullity of the Jacobi operator defined above.

\begin{theorem}\label{thm:main}
For the complex and quaternionic Hopf fibrations with $p,q\ge0$, and for the octonionic Hopf fibration with $p=q=0$, $\Sigma_{p,q}$ is minimal in $(S^{m+1},g_t)$ for every $t>0$. Its unit normal is horizontal, and
\[
|A_t|^2=m+2r(t-1),\qquad
\operatorname{Ric}_{g_t}(\nu,\nu)=m+2r(1-t).
\]
Thus the Jacobi potential is $2m$. Set $t_*=r/(m+r)$. For $0<t\le1$,
\[
\operatorname{ind}(\Sigma_{p,q})=
\begin{cases}1,&0<t\le t_*,\\m+3,&t_*<t\le1.\end{cases}
\]
For $0<t<1$, $t\ne t_*$, the nullity is
\begin{equation}\label{eq:genericnull}
 \operatorname{nul}(\Sigma_{p,q})=d(p+1)(q+1),
\end{equation}
while at the transition and at the round metric, respectively,
\begin{equation}\label{eq:jumpnull}
 \operatorname{nul}_{t_*}(\Sigma_{p,q})=d(p+1)(q+1)+m+2,
\end{equation}
\begin{equation}\label{eq:roundnull}
 \operatorname{nul}_{1}(\Sigma_{p,q})=d^2(p+1)(q+1).
\end{equation}
\end{theorem}

In the octonionic case, $p=q=0$ and $\Sigma_{0,0}=S^7(1/\sqrt2)\times S^7(1/\sqrt2)$. Write $\Sigma=\Sigma_{0,0}$. The next theorem supplements Theorem~\ref{thm:main} by determining the first strictly positive Jacobi eigenvalue and its multiplicity.

\begin{theorem}\label{thm:introoct}
For $0<t\le1$, the smallest strictly positive eigenvalue $\gamma(t)$ of $J_t$ on $\Sigma$ is
\[
\gamma(t)=
\begin{cases}
36,&0<t\le7/57,\\
7/t-21,&7/57<t<1/3,\\
32,&t=1/3,\\
16(t^{-1}-1),&1/3<t<1,\\
4,&t=1.
\end{cases}
\]
Its multiplicity is $35$, $51$, $16$, $168$, $56$, and $70$, respectively, on $(0,7/57)$, $\{7/57\}$, $(7/57,1/3)$, $\{1/3\}$, $(1/3,1)$, and $\{1\}$.
\end{theorem}

The octonionic case also exhibits bifurcation outside the shortening range. The following result is proved in Sections~5--7.

\begin{theorem}\label{thm:introbif}
For the octonionic Hopf squashing at $t=27$, the balanced Clifford hypersurface $\Sigma$ admits a smooth local branch of nonhomogeneous minimal hypersurfaces $\widetilde\Sigma_\varepsilon$. For all sufficiently small nonzero $\varepsilon$,
\[
\operatorname{ind}(\widetilde\Sigma_\varepsilon)=1432,\qquad
\operatorname{nul}(\widetilde\Sigma_\varepsilon)=15.
\]
The branch and its parameter expansion are specified in Theorem~\ref{bif:main}.
\end{theorem}

Al\'ias and Piccione~\cite{AliasPiccione} established
bifurcation and local rigidity results for
constant-mean-curvature Clifford tori in round spheres.
Bifurcation and local rigidity of constant-mean-curvature Clifford tori in three-dimensional Berger spheres were studied by de Lima, de Lira and Piccione~\cite{deLimaDeLiraPiccione}. The present work treats higher-dimensional quaternionic and octonionic Hopf deformations; its spectral calculations and the octonionic minimal-hypersurface branch at $t=27$, including the Morse index and nullity, are distinct from those earlier results.

In the complex case, the index transition was obtained by Torralbo and Urbano~\cite{TorralboUrbano}; their fiber-length parameter $\tau$ corresponds to $t=\tau^2$. The underlying deformation of Riemannian submersion metrics is treated in~\cite{Besse}, and the octonionic Hopf fibration and its symmetries are discussed in~\cite{OPPV}.

Section~2 computes the Jacobi potential, and Section~3 determines the index, nullity, and low spectral data. Section~4 determines the Morse index and nullity at the distinguished Einstein parameters of the quaternionic and octonionic Hopf deformations. Sections~5--6 establish the octonionic bifurcation and its index and nullity; Section~7 contains the required moment calculations.

\subsection*{Acknowledgements}
ChatGPT (OpenAI) was used to assist with language editing, literature searches, and the verification of mathematical calculations. These things were independently reviewed by the author, who takes full responsibility for the content of the manuscript.

\section{Hopf squashing and the Jacobi potential}
We compute the second fundamental form and normal Ricci curvature of the balanced Clifford hypersurfaces under Hopf squashing. Their contributions to the Jacobi potential cancel, so the potential is independent of the squashing parameter.

\begin{proposition}\label{prop:hopfclifford}
Consider one of the three normed-division-algebra Hopf fibrations of a unit round sphere, with fiber dimension \(r\), and let \(\Sigma_{p,q}\) be the associated balanced Clifford hypersurface of dimension \(m\).  Under Hopf squashing \(g_t=g_{\mathcal H}+t g_{\mathcal V}\),
\[
 |A_t|^2-|A_1|^2=2r(t-1),\qquad
 \Ric_{g_t}(\nu,\nu)-\Ric_{g_1}(\nu,\nu)=-2r(t-1).
\]
In particular, \(\Sigma_{p,q}\) is minimal for every \(t>0\) and its Jacobi potential equals \(2m\).
\end{proposition}

\begin{proof}
Write a point of the ambient unit sphere as $(x,y)$ according to the splitting
\[
 \F^N=\F^{p+1}\oplus\F^{q+1}.
\]
Set $c=\cos\theta$ and $s=\sin\theta$ and consider
\[
 \Sigma_\theta=S^a(c)\times S^b(s).
\]
Write $x=cu$ and $y=sv$, where $u\in S^a(1)$ and $v\in S^b(1)$.
The round unit normal to $\Sigma_\theta$ is then
\[
 \nu=(-s u,c v).
\]
For all three Hopf fibrations this normal is horizontal.  In the complex and quaternionic cases this follows from diagonal right multiplication by unit scalars.  In the octonionic case the assertion follows at $(c,s)\in\mathbb O^2$ from the description of the vertical space below, and homogeneity extends it to the orbit.  Hence $|\nu|_{g_t}=1$ for every $t$.

Let $r=d-1$ be the fiber dimension.  In the complex and quaternionic cases, let $I$ be a unit imaginary scalar and use the tangent vectors
\[
 E_1=(cuI,0),\qquad E_2=(0,svI).
\]
The vector $E_1+E_2$ is unit and vertical (for unit imaginary $I$).  Since Hopf squashing multiplies the metric by $t$ in the vertical direction and leaves the horizontal metric unchanged, the restriction of $g_t$ to $\operatorname{span}\{E_1,E_2\}$ is
\begin{equation}\label{eq:block}
 G_\theta=
 \begin{pmatrix}
 c^2+(t-1)c^4&(t-1)c^2s^2\\
 (t-1)c^2s^2&s^2+(t-1)s^4
 \end{pmatrix}.
\end{equation}
Thus
\[
\det G_\theta=t c^2s^2.
\]
There are $r$ such blocks.  The remaining $a-r$ directions in the first factor and $b-r$ directions in the second are horizontal.  Therefore
\[
d\operatorname{vol}_{\Sigma_\theta,g_t}
 =C t^{r/2}c^a s^b\,d\mu,
\]
where $C$ and $d\mu$ are independent of $\theta$.  Since $\partial_\theta=\nu$, the minimal orbit is characterized by
\[
 -a\tan\theta+b\cot\theta=0,
\]
that is,
\begin{equation}\label{eq:balancedradii}
 c^2=\frac am,\qquad s^2=\frac bm.
\end{equation}
Thus the value of $\theta$ for which $\Sigma_\theta$ is minimal is independent of $t$. Equivalently, the mean curvature is $-a\tan\theta+b\cot\theta$, so $\operatorname{tr}A_t=0$ at the balanced radii.

In the complex and quaternionic cases the scalar Hopf action gives the block globally.  In the octonionic case there is no global $S^7$ action; instead use the Hopf map
\[
 h(x,y)=\bigl(|x|^2-|y|^2,\,2x\bar y\bigr).
\]
At the real point $(c,s)\in S^{15}\subset\mathbb O^2$, its vertical tangent space is
\[
\mathcal V_{(c,s)}=\{(cu,su):u\in\operatorname{Im}\mathbb O\}.
\]
For $u\in\operatorname{Im}\mathbb O$,
\[
 dh_{(c,s)}(cu,su)
 =\bigl(0,2cs(u+\bar u)\bigr)=0,
\]
and these seven independent vectors exhaust the kernel because the Hopf fiber has dimension seven.  Consequently the vertical projection of the two tangent vectors $(cu,0)$ and $(0,su)$ is read off from the single vertical vector $(cu,su)$, and Hopf squashing gives the block \eqref{eq:block}.

The hypersurface $S^7(c)\times S^7(s)$ is a $\operatorname{Spin}(8)$ orbit with stabilizer $G_2$, so the pointwise computation at $(c,s)$ determines the homogeneous geometry everywhere.  At $c=s=1/\sqrt2$ the restricted Hopf map is the homogeneous fibration
\[
 S^7=\operatorname{Spin}(7)/G_2
 \longrightarrow
 \operatorname{Spin}(8)/G_2
 \longrightarrow
 \operatorname{Spin}(8)/\operatorname{Spin}(7)=S^7.
\]
The identifications $\operatorname{Spin}(7)/G_2\cong S^7$ and
$\operatorname{Spin}(8)/G_2\cong S^7\times S^7$ are standard consequences of triality; see the Clifford-algebra background in~\cite{Porteous}.  With the convention $A_t(U)=-\nabla^{g_t}_U\nu$, the $a-r$ ordinary directions in the first factor have principal curvature $s/c$, and the $b-r$ ordinary directions in the second have principal curvature $-c/s$. Since $\partial_\theta=\nu$ is a unit normal orthogonal to the hypersurface orbits, differentiating the induced metric in the normal direction gives, for each Hopf block,
\[
 A_{\rm block}=-\frac12G_\theta^{-1}G_\theta'
\]
gives
\[
\operatorname{tr}(A_{\rm block}^2)
 =\frac1{c^2s^2}+2t-4.
\]
Hence
\begin{align*}
 |A_t|^2
 &=(a-r)\frac{s^2}{c^2}+(b-r)\frac{c^2}{s^2}
 +r\left(\frac1{c^2s^2}+2t-4\right).
\end{align*}
Substituting \eqref{eq:balancedradii} and $m=a+b$ gives
\begin{equation}\label{eq:shapeuniv}
 |A_t|^2=m+2r(t-1).
\end{equation}
At $t=1$ this is $|A|^2=m$, the classical value for a minimal Clifford hypersurface in the unit sphere.

For the normal Ricci term, let $X$ be a unit horizontal vector.  Let $\mathcal A$ denote O'Neill's integrability tensor of the Hopf submersion.  For horizontal vector fields $U,V$, it is defined by $\mathcal A_UV=\frac12[U,V]^{\mathcal V}$, where the superscript $\mathcal V$ denotes orthogonal projection onto the vertical distribution.  For the complex, quaternionic, and octonionic Hopf fibrations it satisfies
\begin{equation}\label{eq:Asum}
 \sum_j |\mathcal A_XE_j|^2=r
\end{equation}
where $(E_j)$ is an orthonormal basis of $X^\perp\cap\mathcal H$.  In the complex and quaternionic cases, if $J_\alpha$, $1\le\alpha\le r$, denote the standard complex or quaternionic structures determined by the Hopf fibration, then the vertical components of $\mathcal A_XY$ are $\langle J_\alpha X,Y\rangle$.  Hence
\[
 \sum_j|\mathcal A_XE_j|^2
 =\sum_{\alpha=1}^r |J_\alpha X|^2=r;
\]
see also \cite{Besse,Pipoli}.  In the octonionic case \eqref{eq:Asum} follows from the normalization of the Hopf submersion: the total sphere has sectional curvature $1$ and the base is the round $S^8$ of sectional curvature $4$ \cite{Besse,BKN}.  O'Neill's horizontal curvature formula therefore gives
\[
 4=1+3|\mathcal A_XY|^2
\]
for orthonormal horizontal $X,Y$.  Thus $|\mathcal A_XY|=|Y|$ on $X^\perp\cap\mathcal H$.  Polarization shows that $Y\mapsto\mathcal A_XY$ preserves inner products.  Since $X^\perp\cap\mathcal H$ and $\mathcal V$ are both seven-dimensional, it is an isometry, proving \eqref{eq:Asum} for $r=7$.

The sectional-curvature formulas for Hopf squashing with totally geodesic fibers give
\[
 K_t(X,E_j)=1+3(1-t)|\mathcal A_XE_j|^2,
 \qquad
 K_t(X,t^{-1/2}U_\alpha)=t.
\]
There are $m-r$ horizontal directions transverse to $X$ and $r$ vertical directions.  Summing and using \eqref{eq:Asum} gives
\begin{equation}\label{eq:Ricuniv}
 \Ric_{g_t}(X,X)
 =m-r+3r(1-t)+rt
 =m+2r(1-t).
\end{equation}
The normal $\nu$ is horizontal, so \eqref{eq:shapeuniv} and \eqref{eq:Ricuniv} prove
\[
 |A_t|^2+\Ric_{g_t}(\nu,\nu)=2m.
\]
\end{proof}

\section{The stability transition}
We determine the index and nullity by comparing the eigenvalues of the induced Laplacian with the constant Jacobi potential $2m$ found in Proposition~\ref{prop:hopfclifford}. Starting from the spectrum of the round product at $t=1$, we examine how these eigenvalues change when the Hopf fibers are shortened ($0<t<1$). For the octonionic hypersurface, we also identify the first positive Jacobi eigenvalue and its multiplicity.

\subsection{The nonpositive Jacobi spectrum}
We use the nonnegative Laplacian $\Delta=-\operatorname{div}\nabla$. By Proposition~\ref{prop:hopfclifford}, the Jacobi operator is $J_t=\Delta_t-2m$.

At $t=1$ the induced metric is the round product
\[
 S^a\!\left(\sqrt{a/m}\right)\times
 S^b\!\left(\sqrt{b/m}\right).
\]
Its eigenvalues are
\begin{equation}\label{eq:roundspec}
 \lambda_{k,l}^{(1)}
 =\frac ma k(k+a-1)+\frac mb l(l+b-1).
\end{equation}
The only levels at or below $2m$ are
\[
 (0,0):0,\qquad (1,0),(0,1):m,\qquad (1,1):2m.
\]
Indeed the next single-factor levels already satisfy
\[
 \lambda_{2,0}^{(1)}=2m\frac{a+1}{a}>2m,
 \qquad
 \lambda_{0,2}^{(1)}=2m\frac{b+1}{b}>2m.
\]

The eigenspaces of bidegrees $(1,0)$ and $(0,1)$ consist of the coordinate functions on the respective sphere factors.

The Hopf fibers contained in $\Sigma_{p,q}$ are totally geodesic.  Let \(\Delta_{\mathcal V}\) denote the nonnegative Laplacian along the fibers.  Hopf squashing gives on functions
\begin{equation}\label{eq:lapvar}
 \Delta_t=\Delta_1+\left(\frac1t-1\right)\Delta_{\mathcal V}.
\end{equation}
In the complex and quaternionic cases the diagonal scalar action is isometric, and $\Delta_{\mathcal V}$ is its fiberwise Laplacian (equivalently, its Casimir operator); hence it commutes with $\Delta_1$.  On either first-harmonic family its eigenvalue is the first eigenvalue of the unit round fiber $S^r$, namely
\[
\Delta_{\mathcal V}=r.
\]
In the octonionic case, the vertical Laplacian is not described by a diagonal scalar action. The following lemma determines its action on the low-degree eigenspaces needed for the index and nullity calculation.

\begin{lemma}\label{lem:octspectral}
For the restricted octonionic fibration
\[
 \operatorname{Spin}(7)/G_2\longrightarrow
 \operatorname{Spin}(8)/G_2\longrightarrow
 \operatorname{Spin}(8)/\operatorname{Spin}(7),
\]
the vertical Laplacian preserves the relevant $\operatorname{Spin}(8)$-modules.  On each of the two $8$-dimensional first-harmonic spaces of
$S^7(1/\sqrt2)\times S^7(1/\sqrt2)$ it has eigenvalue $7$.  Inside the $64$-dimensional mixed $(1,1)$ space, its kernel has dimension exactly $8$.
\end{lemma}

\begin{proof}
Set $G=\operatorname{Spin}(8)$, $H=G_2$, and let $K=\operatorname{Spin}(7)$ be the stabilizer of a unit vector in the vector representation $8_v$.  Then
$G/H\cong S^7\times S^7$, while $G/K\cong S^7$ is the equator in the octonionic Hopf base.  The vertical distribution of $G/H\to G/K$ and the round metric are $G$-invariant.  Hence $\Delta_{\mathcal V}$ is a self-adjoint $G$-equivariant operator.

Denote the vector and two half-spin representations of $\operatorname{Spin}(8)$ by $8_v$, $8_s$, and $8_c$, with $K$ the stabilizer of a unit vector in $8_v$.  Thus
\[
 8_v|_K=\mathbf1\oplus\mathbf7,
\]
whereas both half-spin modules restrict to the same irreducible $8$-dimensional spin representation of $K$:
\[
 8_s|_K\cong 8_c|_K\cong\mathbf8.
\]
This is the standard triality realization of
$\operatorname{Spin}(7)/G_2\cong S^7$; see~\cite{Porteous} for the Clifford-algebra description.  Each half-spin module contains a one-dimensional $H$-fixed subspace, so by Frobenius reciprocity it occurs with multiplicity one in $L^2(G/H)$.  Thus $\Delta_{\mathcal V}$ acts by a scalar on each of the two first-harmonic modules.  The $H$-fixed vector restricts along $K/H$ to the degree-one spherical representation of the unit round $S^7$.  Hence that scalar is the first spherical eigenvalue, namely $7$.

For the mixed space, triality gives
\[
 8_s\otimes 8_c = 8_v\oplus 56_v.
\]
The tensor-square decomposition of the eight-dimensional spin representation of $K$ is
\[
 \mathbf8\otimes\mathbf8
 =\mathbf1\oplus\mathbf7\oplus\mathbf{21}\oplus\mathbf{35},
\]
where $\mathbf{35}=\operatorname{Sym}^2_0(\mathbf8)$ and
$\Lambda^2\mathbf8=\mathbf7\oplus\mathbf{21}$ (see~\cite{Porteous}).
Since $8_v|_K=\mathbf1\oplus\mathbf7$, subtraction from the restricted tensor product gives
\[
 56_v|_K=\mathbf{21}\oplus\mathbf{35}.
\]
We can determine the required $G_2$-fixed dimensions directly. The stabilizer $H=G_2$ of a unit spinor in $\mathbf8$ fixes a line, and its orthogonal complement is the standard seven-dimensional $G_2$-module $W$; thus
$\mathbf8|_H=\mathbf1\oplus W$.
Using the standard decompositions
$\operatorname{Sym}^2 W=\mathbf1\oplus\mathbf{27}$ and
$\Lambda^2 W=\mathbf7\oplus\mathbf{14}$, we obtain
\[
 \begin{aligned}
 \mathbf{35}|_H
 &=\operatorname{Sym}^2_0(\mathbf1\oplus W)
   =\mathbf1\oplus\mathbf7\oplus\mathbf{27},\\
 \mathbf{21}|_H
 &=\Lambda^2(\mathbf1\oplus W)-\mathbf7
   =\mathbf7\oplus\mathbf{14}.
 \end{aligned}
\]
Here subtraction denotes removal of a direct summand in the representation ring.
Hence $\dim(56_v)^H=1$, while $\dim(8_v)^H=1$ as well.  By Frobenius reciprocity, both representations therefore occur with multiplicity one in $L^2(G/H)$, so $\Delta_{\mathcal V}$ preserves each summand and acts on it by a scalar.  Moreover $56_v$ contains no $K$-fixed vector, whereas $8_v$ contains exactly one.  A function on $G/H$ is constant on the fibers of $G/H\to G/K$ exactly when it is the pullback of a function on $G/K$.  Equivalently, a $G$-module occurring in $L^2(G/H)$ contributes to the vertical kernel exactly when it contains a $K$-fixed vector.  Therefore the scalar on $56_v$ is positive, while it is zero on $8_v$.  The $8_v$ summand is the eight-dimensional space generated by the coordinate functions of the restricted octonionic Hopf map.  Consequently the fiber-constant part of the mixed eigenspace is precisely the $8_v$ summand and has dimension $8$.
\end{proof}

In all three cases the common first-harmonic eigenvalue is
\begin{equation}\label{eq:first}
 \lambda_{\mathrm{lin}}(t)=m+r\left(\frac1t-1\right).
\end{equation}
It reaches the Jacobi threshold $2m$ precisely when
\[
 m+r\left(\frac1t-1\right)=2m,
\]
that is,
\[
t=t_*:=\frac{r}{m+r}.
\]

For $0<t\le1$, equation \eqref{eq:lapvar} implies $\Delta_t\ge\Delta_1$ as quadratic forms.  In the complex and quaternionic cases $\Delta_{\mathcal V}$ commutes with $\Delta_1$.  In the octonionic case Lemma~\ref{lem:octspectral} shows that the round eigenspaces with eigenvalues $m$ and $2m$ split into $\Delta_{\mathcal V}$-invariant irreducible summands; constants are invariant as well.  Since $\Delta_{\mathcal V}$ is self-adjoint, the orthogonal complement of these low round eigenspaces is also invariant.  On that complement $\Delta_1>2m$, and hence $\Delta_t>2m$.  Thus only the round levels listed above can contribute to the nonpositive Jacobi spectrum.  Constants are the only unstable modes for $0<t\le t_*$, while for $t_*<t\le1$ the two first-harmonic families also become unstable.  Their combined dimension is
\[
 (a+1)+(b+1)=dN=m+2.
\]
Hence
\[
 \ind=1\quad(0<t\le t_*),
\qquad
 \ind=1+(m+2)=m+3\quad(t_*<t\le1).
\]
At $t=t_*$ those same $m+2$ modes enter the kernel rather than the negative spectrum.

It remains to determine the nullity. For $0<t<1$, the persistent Jacobi fields are precisely the fiber-constant functions in the round $(1,1)$ eigenspace. By \eqref{eq:lapvar}, for $0<t<1$ a $(1,1)$ mode remains at $2m$ exactly when it is constant along the Hopf fibers.  In the complex and quaternionic cases the mixed $(1,1)$ space consists of real bilinear forms between the two coordinate blocks.  To compute the invariant part, fix one coordinate in each block and write a real bilinear form as
$B(x,y)=\langle Tx,y\rangle_{\mathbb R}$ for a real-linear endomorphism $T$ of $\F$.
For $u$ a unit scalar, invariance under $(x,y)\mapsto(xu,yu)$ is equivalent to
$T(xu)=(Tx)u$ for all $x,u$.  By real linearity, this is equivalent to
commutation with every right multiplication.  Evaluating at $x=1$ then gives
$T(x)=a x$ for a unique $a\in\F$ (associativity is used here, and is valid for
$\F=\mathbb C,\mathbb H$).  Conversely, left multiplication by any $a\in\F$
commutes with right multiplication, so the space of invariant bilinear forms
has real dimension $d$.  Equivalently, it is spanned by the real components
of $x\overline y$, which is invariant since $(xu)\overline{(yu)}
=x(u\overline u)\overline y=x\overline y$.  Thus the invariant part of the mixed eigenspace has real dimension
\[
 d(p+1)(q+1).
\]
In the octonionic case Lemma~\ref{lem:octspectral} shows that the fiber-constant subspace is exactly the $8_v$ summand generated by the equatorial coordinate functions of the octonionic Hopf map.  Its dimension is $8=d(p+1)(q+1)$.  This gives \eqref{eq:genericnull}.  Adding the two first-harmonic families at $t_*$ gives \eqref{eq:jumpnull}.  At $t=1$ the entire $(1,1)$ eigenspace returns to the kernel, of dimension
\[
 (a+1)(b+1)=d^2(p+1)(q+1),
\]
which gives \eqref{eq:roundnull}.

\begin{proposition}\label{prop:lowspectrum}
Let \(J_t=\Delta_t-2m\).  For \(0<t<1\), the nonpositive spectrum of \(J_t\) consists of the constant eigenvalue
\[
 -2m
\]
with multiplicity \(1\), the persistent zero eigenspace of dimension
\[
 d(p+1)(q+1),
\]
and, only when \(t\ge t_*\), the first-harmonic eigenvalue
\[
 \mu_1(t)
 =-m+r\left(\frac1t-1\right)
\]
with multiplicity \(m+2\).  Here \(\mu_1(t_*)=0\) and \(\mu_1(t)<0\) for \(t>t_*\).  No other Jacobi eigenvalue is nonpositive.

At \(t=1\), the negative spectrum consists of \(-2m\) with multiplicity \(1\) and \(-m\) with multiplicity \(m+2\), while the zero eigenspace has dimension \(d^2(p+1)(q+1)\).
\end{proposition}

\begin{proof}
Let $E_{\le2m}$ be the sum of the round eigenspaces with eigenvalue at most $2m$.  The decomposition of the round eigenspaces shows that $E_{\le2m}$ is invariant under $\Delta_{\mathcal V}$. Since $\Delta_{\mathcal V}$ is self-adjoint, its orthogonal complement is invariant as well.  On $E_{\le2m}^{\perp}$ one has $\Delta_1>2m$, while \eqref{eq:lapvar} gives $\Delta_t\ge\Delta_1$ for $0<t\le1$.  Hence $\Delta_t>2m$ on $E_{\le2m}^{\perp}$.  It remains only to inspect $E_{\le2m}$.  Equation \eqref{eq:first} treats the two first-harmonic spaces.  On the mixed $(1,1)$ space, the kernel of $\Delta_{\mathcal V}$ gives the zero eigenspace of $J_t$, while every positive vertical eigenvalue raises the corresponding Laplace eigenvalue strictly above $2m$ when $t<1$.  The listed eigenvalues therefore exhaust the nonpositive Jacobi spectrum.
\end{proof}

\subsection{The first positive eigenvalue in the octonionic case}
To determine the first positive Jacobi eigenvalue, we analyze the vertical
Laplacian on the low spherical bidegrees and control the remaining spectrum
by a min--max estimate.

\subsubsection{Vertical projections and the spectral cutoff}\label{sec:octonionic-projections}
The following projections establish the invariance required for the min--max argument, without assuming invariance of every spherical bidegree.

Write the balanced hypersurface as pairs $(X,Y)$ of unit octonions, with the
ambient coordinates $(X,Y)/\sqrt2$.  The restricted Hopf map has quotient
coordinate $z=X\bar Y$, and its fiber over $z$ is
\[
 F_z(v)=(zv,v),\qquad v\in S^7.
\]
Here and below products of three octonions are parenthesized as displayed.
The round product Laplacian has eigenvalue
\[
 \lambda_{k,l}=2k(k+6)+2l(l+6)
\]
on $E_{k,l}=\mathcal H_k(S^7)\otimes\mathcal H_l(S^7)$.
The vertical Laplacian has fiber eigenvalues $j(j+6)$.

\begin{lemma}\label{lem:lowoct}
Each $E_{k,l}$ with $\lambda_{k,l}<64$ is invariant under the
vertical Laplacian. In particular, on $E_{2,1}$ its only vertical
eigenvalues are $7$ and $27$. The sum
$L=\bigoplus_{\lambda_{k,l}<64}E_{k,l}$ and its orthogonal complement
are invariant under $\Delta_{\mathcal V}$ and $J_t$.
\end{lemma}
\begin{proof}
On each fiber, a harmonic polynomial in $X$ alone of degree $k$ remains
harmonic of degree $k$, since left multiplication by $z$ is orthogonal.
This proves the assertion for $(k,l)=(k,0),(0,k)$ with $k\le3$.
For $(1,1)$, the fiberwise constant projection of
$f(X,Y)=\langle a,X\rangle\langle b,Y\rangle$ is
\[
 (P_0f)(X,Y)=\frac18\langle a,(X\bar Y)b\rangle.
\]
It is bilinear, hence belongs to $E_{1,1}$.  The complement has vertical
degree two.

For $(2,1)$ take $f(X,Y)=\langle SX,X\rangle\langle w,Y\rangle$,
where $S$ is symmetric and trace-free.  With $A_z=L_{\bar z}SL_z$,
the fiber restriction is $\langle A_zv,v\rangle\langle w,v\rangle$.
The harmonic decomposition of this cubic has degree-one component
$\frac15\langle A_zw,v\rangle$, since
\[
 \Delta_{\mathbb R^8}\bigl(\langle A_zv,v\rangle\langle w,v\rangle\bigr)
 =4\langle A_zw,v\rangle,
 \qquad \Delta_{\mathbb R^8}(|v|^2\ell(v))=20\ell(v)
\]
for a linear function $\ell$.  Consequently
\begin{equation}\label{eq:explicitP7}
 (P_7f)(X,Y)=\frac15\langle S((X\bar Y)w),X\rangle.
\end{equation}
This is linear in $Y$.  To check harmonicity in $X$, set
$T_{Y,w}=R_wR_{\bar Y}$, where $R_a(X)=Xa$.
The composition-algebra polarization identity gives
\[
 T_{Y,w}+T_{Y,w}^*=2\langle w,Y\rangle I.
\]
Since $S$ is symmetric and trace-free,
\[
 \Delta_X\langle ST_{Y,w}X,X\rangle
 =2\operatorname{tr}(ST_{Y,w})=0.
\]
Thus $P_7f\in E_{2,1}$.  The complementary fiberwise cubic harmonic
has degree three.  The $(1,2)$ case follows by interchanging the factors
and using the corresponding right/left multiplication identities.
These are all the bidegrees below $64$.
Since $\Delta_{\mathcal V}$ is self-adjoint, the invariance of $L$
implies the invariance of $L^\perp$. The round Laplacian $\Delta_1$
preserves both spaces by their spectral definition, and hence so does
$J_t=\Delta_1+(t^{-1}-1)\Delta_{\mathcal V}-28$.
No invariance of arbitrary round bidegrees is asserted.
\end{proof}

\begin{lemma}\label{lem:oct-ranks}
In the balanced octonionic case, the vertical eigenspace dimensions
inside the low round bidegrees are
\[
\begin{array}{c|c|c}
(k,l)&\text{vertical eigenvalue}&\text{dimension}\\ \hline
(0,0)&0&1\\
(1,0),(0,1)&7&8+8\\
(1,1)&0,16&8,56\\
(2,0),(0,2)&16&35+35\\
(2,1),(1,2)&7,27&56+56,224+224\\
(3,0),(0,3)&27&112+112
\end{array}
\]
\end{lemma}
\begin{proof}
The dimensions of spherical harmonics on $S^7$ of degrees $0,1,2,3$
are $1,8,35,112$.  The constant-fiber part of $E_{1,1}$ is the
space of linear harmonics in $z=X\bar Y$, hence has dimension eight.
The dimensions for $E_{2,1}$ and $E_{1,2}$ follow from the
independent trace computation in Proposition~\ref{prop:octcollision}.
\end{proof}

The preceding decomposition identifies the low-degree eigenvalues.
The next proposition bounds the remaining spectrum below by $36$.

\begin{proposition}\label{prop:oct-minmax}
For $0<t\le1$, every Jacobi eigenvalue arising from the orthogonal
complement of $L=\bigoplus_{\lambda_{k,l}<64}E_{k,l}$ is at least $36$.
For $0<t<1$, the eigenspace at $36$ in this complement has dimension $35$.
\end{proposition}
\begin{proof}
Lemma~\ref{lem:lowoct} makes $L$ invariant under the vertical
Laplacian, hence under $J_t=\Delta_1+(t^{-1}-1)\Delta_{\mathcal V}-28$.
Self-adjointness gives invariance of $L^\perp$.  Since
$t^{-1}-1\ge0$ and $\Delta_1\ge64$ on $L^\perp$, we have
$J_t\ge36$ there.
For $t<1$, equality requires both round eigenvalue $64$ and vertical
eigenvalue zero.  The only round bidegree at $64$ is $E_{2,2}$.
The vertically constant functions are pullbacks of functions of
$z=X\bar Y$.  A quotient harmonic quadratic $Q(z)$ pulls back to
$E_{2,2}$: it is quadratic and harmonic separately in $X$ and $Y$,
because the multiplication maps by a fixed unit octonion are orthogonal.
Conversely a vertically constant function in the round $64$-eigenspace
is a quotient eigenfunction of degree two, since the quotient is the
round $S^7(1/2)$.  Its dimension is $35$.
\end{proof}

\begin{proposition}
\label{prop:octcollision}
For the balanced octonionic Clifford hypersurface at $t=1/3$, the
Jacobi eigenvalue $32$ has multiplicity $168$.  Its eigenspace is the
orthogonal sum of the vertical-degree-two component of $E_{1,1}$ and
the vertical-degree-one components of $E_{2,1}$ and $E_{1,2}$.
\end{proposition}

\begin{proof}
Here $E_{k,\ell}$ denotes the product of spherical harmonic spaces of
degrees $k$ and $\ell$ on the two spheres of radius $1/\sqrt2$.  The
round eigenvalue is
\[
 \lambda_{k,\ell}=2k(k+6)+2\ell(\ell+6).
\]
For $a\in S^7\subset\mathbb O$, the Hopf fiber over $a$ is
\[
 F_a=\{(au,u)/\sqrt2:u\in S^7\}.
\]
This parametrization is an isometry from the unit round $S^7$.
Consequently the restriction of an element of $E_{2,1}$ to a fiber is
a polynomial of degree three and decomposes into spherical harmonics
of degrees one and three.  The only possible vertical Laplace
eigenvalues on $E_{2,1}$ are therefore $7$ and $27$.

Their multiplicities can be computed without assuming that the fibers are
orbits of a group action.  For the degree-$k$ harmonic space $\mathcal
H_k$ on a round $S^7(R)$, with an orthonormal basis $\{\phi_\alpha\}$,
the addition formula and isotropy give, for every tangent vector $v$,
\[
 \sum_\alpha |d\phi_\alpha(v)|^2
 =\frac{\dim\mathcal H_k}{\operatorname{Vol}(S^7(R))}
 \frac{k(k+6)}{7R^2}|v|^2.
\]
For a unit vertical tangent vector $V=(V_x,V_y)$ along the balanced
Clifford hypersurface, the fiber parametrization yields
$|V_x|^2=|V_y|^2=1/2$.  Applying the preceding identity to a product
orthonormal basis of $E_{2,1}$ and summing over the seven vertical
directions gives
\[
 \operatorname{tr}(\Delta_{\mathcal V}|_{E_{2,1}})
 = (35\cdot8)\left(\frac{32}{2}+\frac{14}{2}\right)
 =6440.
\]
The restrictionn of $\Delta_{\mathcal V}$ to $E_{2,1}$ is well defined:
the required invariance follows from Lemma~\ref{lem:lowoct}.
If $d_7$ and $d_{27}$ denote the multiplicities, then
\[
 d_7+d_{27}=280,\qquad 7d_7+27d_{27}=6440,
\]
so $d_7=56$ and $d_{27}=224$.  The same calculation applies to
$E_{1,2}$.

On $E_{1,1}$, the fiber restrictions are quadratic; the vertical
spectrum is $0$ with multiplicity eight and $16$ with multiplicity
$56$.  Since $J_{1/3}=\Delta_1+2\Delta_{\mathcal V}-28$, the
vertical-degree-two part of $E_{1,1}$ and the vertical-degree-one
parts of $E_{2,1}$ and $E_{1,2}$ all have Jacobi eigenvalue $32$.
Their dimensions sum to $56+56+56=168$.

Finally, any further contribution would require
$\lambda_{k,\ell}\le60$.  Besides the three spaces just considered,
the only possibilities are $(0,0)$, $(1,0)$, $(0,1)$, $(2,0)$,
$(0,2)$, $(3,0)$, and $(0,3)$.  Their fiber restrictions have,
respectively, vertical degrees $0$, $1$, $1$, $2$, $2$, $3$, and $3$.
Substitution into $\lambda_{k,\ell}+2\lambda_{\mathcal V}-28$
shows that none has eigenvalue $32$.
\end{proof}

\begin{theorem}\label{thm:octgap}
For the balanced octonionic Clifford hypersurface and $0<t\le1$, let
$\gamma(t)$ be the smallest strictly positive eigenvalue of
$J_t=\Delta_t-28$. Then
\[
\gamma(t)=
\begin{cases}
36,&0<t\le 7/57,\\
7/t-21,&7/57<t<1/3,\\
32,&t=1/3,\\
16(t^{-1}-1),&1/3<t<1,\\
4,&t=1.
\end{cases}
\]
Its multiplicity is respectively $35$, $51$, $16$, $168$, $56$, and
$70$ on the six sets
$(0,7/57)$, $\{7/57\}$, $(7/57,1/3)$, $\{1/3\}$,
$(1/3,1)$, and $\{1\}$.
\end{theorem}
\begin{proof}
By Lemma~\ref{lem:lowoct}, the low-energy space and its orthogonal
complement reduce $J_t$. Proposition~\ref{prop:oct-minmax} bounds the
spectrum on the complement below by $36$. The vertical eigenvalues
and their multiplicities in the low bidegrees are given by
Lemma~\ref{lem:oct-ranks} and Proposition~\ref{prop:octcollision}.
The complete list of Jacobi branches
below round level $64$, apart from the constant branch $-28$, is
\[
\begin{array}{c|c}
E_{1,0}\oplus E_{0,1}&7/t-21\\
E_{1,1}&0,\quad16/t-16\\
E_{2,0}\oplus E_{0,2}&16/t-12\\
E_{2,1}\oplus E_{1,2}&7/t+11,\quad27/t-9\\
E_{3,0}\oplus E_{0,3}&27/t-1.
\end{array}
\]
The pullbacks of degree-two harmonics from the base $S^7(1/2)$
have round eigenvalue $64$ and vertical eigenvalue zero, hence
Jacobi eigenvalue $36$ with multiplicity $35$. For $t<1/3$,
$7/t-21>0$, and direct comparison of the displayed branches shows
that the first positive eigenvalue is $\min\{36,7/t-21\}$.
The two candidates coincide at $t=7/57$, where their multiplicities
add. At $t=1/3$ the first-harmonic branch becomes zero and the next
positive value is $32$, contributed by three $56$-dimensional
summands as in Proposition~\ref{prop:octcollision}.
For $1/3<t<1$, the first-harmonic branch is negative and
$16/t-16$ is positive. It is smaller than $16/t-12$ by $4$,
smaller than $7/t+11$ because $27-9/t>0$, and smaller than
$36$; all remaining displayed branches are larger. Its multiplicity
is $56$. At $t=1$, this branch joins the kernel; the first positive
round Jacobi eigenvalue is $32-28=4$, on
$E_{2,0}\oplus E_{0,2}$, of dimension $35+35=70$.
\end{proof}

\subsection{Quaternionic representation-theoretic lemmas}\label{sec:quaternionic-lemmas}
We first determine the vertical spectrum of mixed first harmonics and the dimension of the fiber-invariant quadratic harmonics.

\begin{lemma}\label{lem:quatmixed}
For the diagonal right $\mathrm{Sp}(1)$ action on
$\mathbb H^{p+1}\oplus\mathbb H^{q+1}$, the real mixed bidegree
$(1,1)$ has vertical Laplace eigenvalues $0$ and $8$, of respective
multiplicities $4(p+1)(q+1)$ and $12(p+1)(q+1)$.
\end{lemma}
\begin{proof}
It suffices to work with one pair of quaternionic coordinates.  Every
real bilinear form is uniquely of the form
$B(x,y)=\langle Ax,y\rangle$ with $A\in\operatorname{End}_{\mathbb R}(\mathbb H)$.
Under simultaneous right multiplication by a unit quaternion $u$,
$A$ is conjugated by $R_u$.  The commuting endomorphisms are exactly
$L_a$, $a\in\mathbb H$, a four-dimensional real subspace.  Using the
orthogonal decomposition
$\operatorname{End}_{\mathbb R}(\mathbb H)
=\bigoplus_{a\in\{1,i,j,k\}}L_a\operatorname{span}_{\mathbb R}
\{I,R_i,R_j,R_k\}$, the orthogonal complement consists of four copies
of the adjoint three-dimensional representation of $\mathrm{Sp}(1)$.
Under conjugation by unit quaternions, the span of
$R_i,R_j,R_k$ carries the adjoint representation. Its matrix
coefficients on the unit $3$-sphere are spherical harmonics of degree
two and therefore have vertical Laplace eigenvalue $2(2+2)=8$;
the commuting summand has eigenvalue zero.
Taking the direct sum over the $(p+1)(q+1)$ coordinate pairs proves
the assertion.
\end{proof}

\begin{lemma}\label{lem:quatquadratic}
Let $n\ge1$.  On the unit sphere in $\mathbb H^n$, the right
$\mathrm{Sp}(1)$-invariant real quadratic spherical harmonics form a
space of dimension $n(2n-1)-1$.
\end{lemma}
\begin{proof}
Every invariant real quadratic form is $x\mapsto\langle Hx,x\rangle$
for a quaternionic Hermitian $n\times n$ matrix $H$.  Such matrices
have $n$ real diagonal parameters and $4\binom n2$ real off-diagonal
parameters, giving dimension $n+2n(n-1)=n(2n-1)$.  The Euclidean Laplacian of the quadratic form is a nonzero constant
multiple of the real trace of $H$. Thus harmonicity is equivalent to
$\operatorname{tr}_{\mathbb H}H=0$, which imposes one real linear
condition.
\end{proof}

\begin{remark}\label{rem:riccihorizontal}
For Hopf squashing of a Hopf submersion with unit round total
space and totally geodesic unit round fibers, O'Neill's formulas give,
for horizontal orthonormal $X,Y$ and vertical unit $U$ with respect to
$g_1$,
\[
 K_t(X,Y)=1+3(1-t)|\mathcal A_XY|^2,\qquad
 K_t(X,t^{-1/2}U)=t|\mathcal A_XU|^2.
\]
The Hopf tensor satisfies
$\sum_{Y\perp X}|\mathcal A_XY|^2=r$ and
$\sum_U|\mathcal A_XU|^2=r$; the second equality follows by
adjointness of the horizontal--vertical O'Neill maps.  Since the
ambient dimension is $m+1$, the horizontal orthogonal complement of
$X$ has dimension $m-r$, whence
\[
 \Ric_{g_t}(X,X)=(m-r)+3r(1-t)+rt=m+2r(1-t).
\]
Thus the horizontal Ricci curvature depends only on the summed
O'Neill tensor norms, not on their values on individual horizontal
two-planes. In particular, the normal Ricci formula used in
Proposition~\ref{prop:hopfclifford} applies to both the quaternionic
and exceptional octonionic fibrations.
\end{remark}

\subsection{The first positive quaternionic Jacobi eigenvalue}

We determine the least strictly positive Jacobi eigenvalue and its multiplicity
for the quaternionic Hopf deformation. The answer is piecewise, reflecting
changes in the eigenspace that realizes the spectral gap. We use the
nonnegative Laplace convention and restrict to $0<t\le1$.

\begin{theorem}\label{thm:quatgap}
Let $a=4p+3$, $b=4q+3$, $m=a+b$, and assume $p+q\ge1$. Set
\[
 M=\max\{a,b\},\qquad \delta=\frac{2m}{M},\qquad
 t_c=\frac3{m+3+\delta},\quad t_0=\frac3{m+3},\quad
 t_b=\frac8{8+\delta}.
\]
Here $t_c$ and $t_b$ mark crossings with the constant positive branch
$\delta$, while $t_0$ is the zero of the first-harmonic branch.
Then $0<t_c<t_0<t_b<1$, and the least strictly positive
Jacobi eigenvalue is
\[
 \mu_+(t)=\begin{cases}
 \delta,&0<t\le t_c,\\
 3/t-m-3,&t_c<t<t_0,\\
 \delta,&t_0\le t\le t_b,\\
 8(t^{-1}-1),&t_b<t<1,\\
 \delta,&t=1.
 \end{cases}
\]
Let $d_r=r(2r+3)$ and define
\[
 d=\begin{cases}d_p,&p>q,\\d_q,&q>p,\\2d_p,&p=q\ge1.\end{cases}
 \qquad L=m+2,\qquad H=12(p+1)(q+1).
\]
The multiplicities are
\[
\begin{array}{c|c}
\text{Parameter range} & \text{Multiplicity}\\ \hline
0<t<t_c & d\\
t=t_c & d+L\\
t_c<t<t_0 & L\\
t=t_0 & d\\
t_0<t<t_b & d\\
t=t_b & d+H\\
t_b<t<1 & H
\end{array}
\]
At $t=1$ the multiplicity is $D(M)$ if $a\ne b$ and
$2D(M)$ if $a=b$, where $D(M)=M(M+3)/2$.
\end{theorem}

\begin{proof}
The round product Laplacian on the biharmonic space
$E_{k,l}=\mathcal H_k(S^a)\otimes\mathcal H_l(S^b)$ has eigenvalue
\[
 \lambda_{k,l}=\frac ma k(k+a-1)+\frac mb l(l+b-1).
\]
Outside $E_{0,0}$, $E_{1,0}$, $E_{0,1}$, and $E_{1,1}$,
its eigenvalues are at least $2m+\delta$. Equality is possible
only in $E_{2,0}$ when $a=M$ or $E_{0,2}$ when $b=M$.
The diagonal right $\mathrm{Sp}(1)$ action preserves each $E_{k,l}$,
and its Casimir is the vertical Laplacian. Consequently
$J_t=\Delta_1+(t^{-1}-1)\Delta_{\mathcal V}-2m$ is at least
$\delta$ on the sum of all these higher biharmonic spaces.
For $t<1$, equality requires vertical invariance, since
$t^{-1}-1>0$ and the vertical Laplacian is nonnegative. Such invariant
quadratics are precisely the trace-free quaternionic Hermitian forms
on the larger factor, of real dimension
$(r+1)(2r+1)-1=d_r$ for $r=p$ or $q$.

Constants give $-2m$. The two first-harmonic spaces have combined
dimension $(a+1)+(b+1)=L$ and Jacobi eigenvalue
$A(t)=3/t-m-3$. For the mixed space, on each of the
$(p+1)(q+1)$ quaternionic coordinate pairs the diagonal right action
on real bilinear forms decomposes as four trivial representations
and four copies of the three-dimensional adjoint representation.
Thus the vertical eigenvalues are $0$ and $8$, of multiplicities
$4(p+1)(q+1)$ and $H$, respectively. The corresponding Jacobi
eigenvalues are $0$ and $B(t)=8(t^{-1}-1)$.

For $p+q\ge1$ one has $A(t_c)=\delta$, $A(t_0)=0$,
$B(t_b)=\delta$, with $A$ and $B$ strictly decreasing.
Also $B(t)>A(t)$ for $t\le t_0$, since
$B(t)-A(t)=5/t+m-5>0$.
For $t\ge t_0$, $A(t)\le0$; and for $t>t_0$,
$B(t)>0$ until $t=1$. These comparisons yield the stated
piecewise minimum among the strictly positive branches.
The multiplicities follow from the equality cases just described.
At $t=1$ the vertical term vanishes, so all degree-two spherical
harmonics on the larger factor contribute; their dimension is
$\dim\mathcal H_2(S^M)=M(M+3)/2$.
\end{proof}

When $p=q=0$, the fiber-invariant quadratic modes used in the
preceding theorem are absent. The first positive eigenvalue therefore
requires a separate analysis.

\begin{proposition}\label{prop:quat33gap}
For $p=q=0$, the least strictly positive Jacobi eigenvalue is
\[
 \mu_+(t)=\begin{cases}
 20,&0<t\le3/29,\\
 3/t-9,&3/29<t<1/3,\\
 16,&t=1/3,\\
 8(t^{-1}-1),&1/3<t<1,\\
 4,&t=1.
 \end{cases}
\]
Its multiplicities on the successive strata are
$9,17,8,36,12,18$.
\end{proposition}

\begin{proof}
Here $\lambda_{k,l}=2k(k+2)+2l(l+2)$ and
$J_t=\Delta_1+(t^{-1}-1)\Delta_{\mathcal V}-12$.
The quotient of $S^3\times S^3$ by diagonal right multiplication
is the round $S^3(1/2)$. Its degree-two spherical harmonics
pull back to a nine-dimensional vertical kernel with round
eigenvalue $32$, hence Jacobi eigenvalue $20$.
The linear spaces contribute $A(t)=3/t-9$ with multiplicity eight.
The mixed space contributes zero with multiplicity four and
$B(t)=8(t^{-1}-1)$ with multiplicity twelve.
The next potentially competing round spaces have eigenvalues
$16$ (pure quadratics), $22$ (bidegrees $(2,1),(1,2)$),
and $30$ (pure cubics). Their positive Jacobi eigenvalues are
respectively $8/t-4$, $7+3/t$ or $15/t-5$, and $15/t+3$.
The $(2,1)$ and $(1,2)$ vertical degree-one components have
multiplicity twelve each, by the diagonal $\mathrm{Sp}(1)$
Clebsch--Gordan decomposition $V_2\otimes V_1=V_1\oplus V_3$.
Every higher round level is at least $32$, so its Jacobi
eigenvalues are at least $20$ for $t\le1$.
Direct comparison of the displayed branches gives the formula.
At $t=1/3$, the three degree-eight vertical branches in
$E_{1,1}$, $E_{2,1}$, and $E_{1,2}$ each contribute twelve
functions at eigenvalue $16$, giving multiplicity $36$.
At $t=1$ the first positive round Jacobi eigenvalue is $4$,
from the two pure quadratic spaces, each of dimension nine.
\end{proof}

\subsection{Sharpness of the spectral cutoffs}\label{sec:spectral-sharpness}
The spectral cutoff is attained in both cases considered below.

\begin{proposition}
For the balanced octonionic hypersurface and $0<t<1$, the restriction
of $J_t$ to the orthogonal complement of the round eigenspaces below
$64$ has lowest eigenvalue exactly $36$, of multiplicity $35$.
\end{proposition}
\begin{proof}
By Lemma~\ref{lem:lowoct}, the sum of the round eigenspaces below $64$
is invariant under $\Delta_{\mathcal V}$. On its orthogonal complement,
\[
\langle J_t f,f\rangle-36\|f\|^2
=\langle(\Delta_1-64)f,f\rangle
 +(t^{-1}-1)\langle\Delta_{\mathcal V}f,f\rangle\ge0.
\]
For $t<1$, equality holds if and only if $f$ has round eigenvalue
$64$ and is constant on the fibers. The only bidegree with round
eigenvalue $64$ is $(2,2)$. A fiber-constant eigenfunction descends
to the quotient $S^7(1/2)$. The fibers of the restricted Hopf
submersion are totally geodesic and hence minimal; consequently the
Laplacian of a basic function is the pullback of the base Laplacian
(cf.~\cite{Besse}). Thus its eigenvalue agrees with that of the
quotient Laplacian. On $S^7(1/2)$, eigenvalue $64$ corresponds exactly to
harmonic degree two, of dimension $35$. Conversely, the Hopf map is bilinear in the two octonionic
coordinates. The pullback of a quadratic harmonic on the quotient
is therefore homogeneous of degree two in each variable. Since
the pullback is basic, the submersion formula gives round
eigenvalue $64$; the product-sphere eigenvalue on bidegree
$(2,2)$ is also $64$. The separate Laplacian eigenvalues are
therefore those of quadratic harmonics, so the pullback belongs
to $E_{2,2}$. Thus
all $35$ functions occur and there are no others.
\end{proof}

\begin{proposition}
Let $a=4p+3$, $b=4q+3$, $m=a+b$, and $M=\max\{a,b\}$.
Assume $p+q\ge1$ and set $\delta=2m/M$.
For $0<t<1$, the lowest eigenvalue of $J_t$ on the orthogonal
complement of $E_{0,0}\oplus E_{1,0}\oplus E_{0,1}\oplus E_{1,1}$
is $\delta$. Its multiplicity is $d_p$ if $p>q$, $d_q$ if
$q>p$, and $2d_p$ if $p=q$, with $d_r=r(2r+3)$.
\end{proposition}
\begin{proof}
The diagonal right $\mathrm{Sp}(1)$ action preserves each product
bidegree. The first round eigenvalue outside the four indicated
bidegrees is $2m+\delta$, reached precisely by the pure quadratic
bidegrees on the factors of maximal dimension $M$. For $t<1$ we have
\[
\langle(J_t-\delta)f,f\rangle
=\langle(\Delta_1-2m-\delta)f,f\rangle
 +(t^{-1}-1)\langle\Delta_{\mathcal V}f,f\rangle.
\]
Both terms are nonnegative on the indicated complement.
Equality requires separately that $f$ belong to the lowest round
eigenspace on this complement, namely the pure quadratics of a
maximal factor, and that $\Delta_{\mathcal V}f=0$, equivalently
that $f$ be invariant under right $\mathrm{Sp}(1)$. Such quadratic forms are precisely the
trace-free quaternionic Hermitian forms on $\mathbb H^{r+1}$,
where $4r+3=M$. Lemma~\ref{lem:quatquadratic} gives their
dimension $(r+1)(2r+1)-1=r(2r+3)$. If both factors have maximal
dimension, the dimensions add.
\end{proof}

\begin{remark}
The sharpness statements above concern $0<t<1$. At $t=1$ the
vertical contribution vanishes, so the equality eigenspaces may
have larger multiplicity.
\end{remark}

\subsection{The adiabatic limit in the octonionic case}\label{sec:oct-adiabatic}
The low Jacobi spectrum has the following limit as $t\downarrow0$. In the balanced octonionic case the restricted Hopf
map is a Riemannian submersion
\[
 \pi_\Sigma:(\Sigma,g_1)\longrightarrow S^7(1/2)
\]
with totally geodesic unit $S^7$ fibers. Let $P_0$ be fiberwise
averaging, the orthogonal projection onto the basic functions.

\begin{proposition}
\label{prop:oct-adiabatic}
For $0<t\le1$, the basic subspace $P_0L^2(\Sigma)$ and its
orthogonal complement reduce $J_t$. On the basic subspace the spectrum,
including multiplicities, is independent of $t$ and equals
\[
 \bigl\{4j(j+6)-28:j=0,1,2,\ldots\bigr\},
 \qquad
 d_j=\binom{j+7}{7}-\binom{j+5}{7}.
\]
On the orthogonal complement one has the quadratic-form estimate
\[
 J_t\big|_{\ker P_0}\ \ge\ 7(t^{-1}-1)-28.
\]
Consequently, for each fixed $R\in\mathbb R$, all eigenvalues of
$J_t$ below $R$ belong to the basic subspace when $t>0$ is
sufficiently small. The bounded spectrum in the adiabatic limit is
precisely the shifted Laplace spectrum of $S^7(1/2)$.
\end{proposition}
\begin{proof}
Since the fibers are totally geodesic, their mean curvature vanishes.
For every smooth function $u$ on the base, the submersion identity gives
\[
 \Delta_1(u\circ\pi_\Sigma)
   =(\Delta_{S^7(1/2)}u)\circ\pi_\Sigma.
\]
Thus the basic subspace is invariant under $\Delta_1$. Since
$\Delta_1$ is self-adjoint, its orthogonal complement is invariant as
well. The vertical Laplacian vanishes on basic functions, and it
preserves their orthogonal complement by self-adjointness. Since both Laplacians preserve the basic subspace and its orthogonal
complement, this decomposition reduces
$J_t=\Delta_1+(t^{-1}-1)\Delta_{\mathcal V}-28$
for every $0<t\le1$.
On basic functions, $J_t=\Delta_{S^7(1/2)}-28$ under pullback.
The eigenvalues of the sphere of radius $1/2$ are $4j(j+6)$,
with the stated multiplicities.

The vertical fibers are unit round $S^7$ spheres. Their first positive
Laplacian eigenvalue is $7$. Fiberwise Poincar\'e inequality therefore
gives, for $f\perp P_0L^2(\Sigma)$,
\[
 \langle\Delta_{\mathcal V}f,f\rangle\ge7\|f\|^2.
\]
Since $\Delta_1\ge0$ and $t^{-1}-1\ge0$, the asserted bound follows.
Its right-hand side tends to $+\infty$ as $t\downarrow0$.
Consequently, for any fixed bounded spectral interval, all eigenvalues
of $J_t$ in that interval are basic for sufficiently small $t$;
their multiplicities are therefore independent of $t$.
\end{proof}

\begin{corollary}
For sufficiently small $t>0$, the only negative eigenvalue is $-28$
with multiplicity one, the nullity is eight, and the first positive
eigenvalue is $36$ with multiplicity $35$. More precisely, the
spectral formulas above give these conclusions
for $0<t<7/57$.
\end{corollary}
\begin{proof}
The first three basic levels are $-28$, $0$, and $36$, with
multiplicities $1$, $8$, and $35$. Proposition~\ref{prop:oct-adiabatic}
excludes nonbasic eigenvalues in any fixed bounded interval for
sufficiently small $t$. The first-positive-eigenvalue formula in
Theorem~\ref{thm:octgap} sharpens this to the stated parameter range. At $t=7/57$ a nonbasic branch
meets the eigenvalue $36$, so the endpoint is excluded from the
multiplicity statement.
\end{proof}

\section{Quaternionic and octonionic Einstein consequences}
For a $g_t$-unit vertical vector, the Ricci-curvature formula for Hopf squashing gives
\[
\Ric_{g_t}(U,U)=\frac{r-1}{t}+(m+1-r)t;
\]
see \cite{Besse,BKN}.  For the quaternionic Hopf metrics, the transitive group is $\operatorname{Sp}(N)\operatorname{Sp}(1)$ on $S^{4N-1}$, while for the octonionic Hopf metrics it is $\operatorname{Spin}(9)$ on $S^{15}$. In each case, the horizontal and vertical isotropy representations are irreducible and inequivalent. Consequently, the invariant Ricci tensor is scalar on each summand and has no mixed components.  Thus equality of the horizontal and vertical Ricci values is equivalent to the Einstein condition.  Combining the vertical formula with \eqref{eq:Ricuniv}, we obtain
\begin{equation}\label{eq:EinsteinEq}
 m+2r(1-t)=\frac{r-1}{t}+(m+1-r)t.
\end{equation}
Multiplying by $t$ and factoring gives
\[
 (t-1)\bigl((m+r+1)t-(r-1)\bigr)=0,
\]
so the Einstein parameters are $t=1$ and $t=(r-1)/(m+r+1)$.
In the quaternionic case, $r=3$ and $m=4N-2$.  Apart from $t=1$, the solution is Jensen's parameter \cite{Jensen}
\[
 t_J=\frac1{2N+1}<\frac3{4N+1}=t_*.
\]
In the octonionic case, $r=7$ and $m=14$, and \eqref{eq:EinsteinEq} has the roots $1$ and
\[
 t_{\operatorname{Spin}(9)}=\frac3{11}<\frac13=t_*.
\]
The latter is the parameter of the non-round $\operatorname{Spin}(9)$-homogeneous Einstein metric on $S^{15}$ \cite{BK,BKN}.  Both non-round parameters therefore lie in the index-one range of Section~3.

\begin{corollary}\label{cor:Einstein}
Every $S^{4N-1}$ equipped with Jensen's non-round homogeneous Einstein metric carries, for each splitting $N=p+q+2$, the non-totally-geodesic index-one minimal Clifford hypersurface
\[
 S^{4p+3}\!\left(\sqrt{\frac{4p+3}{4N-2}}\right)
 \times
 S^{4q+3}\!\left(\sqrt{\frac{4q+3}{4N-2}}\right).
\]
Its nullity is $4(p+1)(q+1)$.  The non-round $\operatorname{Spin}(9)$-Einstein metric on $S^{15}$ carries the balanced minimal hypersurface
\[
 S^7(1/\sqrt2)\times S^7(1/\sqrt2)
\]
of index one and nullity eight.
\end{corollary}

On $S^{15}$, quaternionic and octonionic Hopf squashing yield distinct non-round homogeneous Einstein metrics, at $t=1/9$ and $t=3/11$, respectively; compare Ziller's classification \cite{Ziller}.  The index and nullity statements used here concern $0<t\leq1$, the range established in Section~3.

\section{Bifurcation at the exceptional parameter \texorpdfstring{$t=27$}{t=27}}
We study birfurcation at $t=27$.
\subsection{The bifurcation theorem and the area functional}
Write $\Sph^{15}\subset\Oct\oplus\Oct$ and consider
\begin{equation}\label{bif:hopf}
\pi(a,b)=(|a|^2-|b|^2,2a\bar b),\qquad g_t=g_{\HH}+t g_{\VV},\quad t>0.
\end{equation}
The balanced Clifford hypersurface is
\begin{equation}\label{bif:sigma}
\Sigma=\{2^{-1/2}(X,Y):X,Y\in\Sph^7\}.
\end{equation}
We use the sign convention that the Jacobi operator $J$ has quadratic form $\int fJf\,dA$, and the Morse index counts negative eigenvalues.
\begin{theorem}\label{bif:main}
There exists a smooth local branch of nontrivial $\operatorname{Spin}(7)$-invariant minimal hypersurfaces $\widetilde\Sigma_\varepsilon$ in $(\Sph^{15},g_{t(\varepsilon)})$, represented as normal graphs over $\Sigma$, with
\begin{equation}\label{bif:expansion}
 u_\varepsilon=\varepsilon p+\varepsilon^2w_2+O(\varepsilon^3),\qquad
 t(\varepsilon)=27-\frac{12746071137}{28028000}\varepsilon^2+O(\varepsilon^4),
\end{equation}
where $p$ and $w_2$ are given in \eqref{bif:pdef} and \eqref{bif:wdef}. For sufficiently small $\varepsilon\ne0$,
\[
\operatorname{ind}(\widetilde\Sigma_\varepsilon)=1432,\qquad
\operatorname{nul}(\widetilde\Sigma_\varepsilon)=15.
\]
\end{theorem}
The proof begins with an exact expression for the area of normal graphs over $\Sigma$.

For $\theta=\pi/4+u$ set $c=\cos\theta$, $s=\sin\theta$, and
\[
F_u(X,Y)=(cX,sY),\qquad z=X\bar Y.
\]
The $\theta$-direction $(-sX,cY)$ is a horizontal unit vector. For tangent vectors $U\in T_X\Sph^7$, $V\in T_Y\Sph^7$, the induced metric on a constant-$\theta$ slice is
\begin{equation}\label{bif:metric}
h_{t,\theta}(U,V)=c^2|U|^2+s^2|V|^2+(t-1)|c^2\bar zU+s^2V|^2.
\end{equation}
After identifying the tangent spaces by $U\mapsto\bar zU$, each of seven blocks has matrix
\[
\begin{pmatrix}c^2+(t-1)c^4&(t-1)c^2s^2\\(t-1)c^2s^2&s^2+(t-1)s^4\end{pmatrix},
\quad \det=tc^2s^2.
\]
Consequently $F_u^*g_t=h_{t,\theta}+du\otimes du$ and
\begin{equation}\label{bif:energy}
\mathcal A_t(u)=C_t\int_{\Sph^7\times\Sph^7}\cos^7(2u)\sqrt{1+Q_t(\theta,\nabla u)}\,dX\,dY,
\quad C_t=2^{-7}t^{7/2},
\end{equation}
where
\begin{align}\label{bif:Q}
Q_t={}&(\tan^2\theta+t^{-1})|\nabla_Xu|^2+(\cot^2\theta+t^{-1})|\nabla_Yu|^2\\
&+2(t^{-1}-1)\av{\bar z\nabla_Xu,\nabla_Yu}.\nonumber
\end{align}
Hereafter brackets $\av{f}$ without two arguments denote integration against normalized product volume; the constant $C_t$ does not affect index or nullity.
\begin{lemma}\label{bif:congruence}
At a minimal graph $u$, the second variation of \eqref{bif:energy} is congruent to the geometric Jacobi quadratic form. In particular, their indices and nullities agree.
\end{lemma}
\begin{proof}
The variation $u+sv$ has normal velocity $\rho_uv$, where
$\rho_u=(1+Q_t(\theta,\nabla u))^{-1/2}>0$. At a minimal immersion tangential components of the velocity do not contribute to the second variation. Thus
\[
D^2\mathcal A_t(u)[v,w]=\int_{\Sigma_u}(\rho_uv)J_u(\rho_uw)\,dA_u.
\]
Since $\rho_u$ is smooth and strictly positive on the compact hypersurface, multiplication by $\rho_u$ is an isomorphism on the relevant Sobolev spaces. The two quadratic forms therefore have the same index and nullity.
\end{proof}
\subsection{The limiting Jacobi spectrum}
Set $L_{\VV}=-\Delta_{\VV}$. The quadratic term of \eqref{bif:energy} gives
\begin{equation}\label{bif:jacobi}
J_t=-\Delta_{\Sigma,1}+(t^{-1}-1)L_{\VV}-28.
\end{equation}
For $\PP_{k,l}=\HH_k(\Sph^7_X)\otimes\HH_l(\Sph^7_Y)$, the round product Laplacian has eigenvalue $2[k(k+6)+l(l+6)]$.
\begin{lemma}\label{bif:bigradeproof}
For every $k,l\ge0$, the vertical Laplacian $L_{\VV}$ preserves
$\PP_{k,l}=\HH_k(\Sph^7_X)\otimes\HH_l(\Sph^7_Y)$.
\end{lemma}
\begin{proof}
Use Euclidean coordinates on $\Oct\oplus\Oct$ and write $E_X=X\cdot\partial_X$,
$E_Y=Y\cdot\partial_Y$. Let
\[
 C_{ij}(X,Y)=\langle(X\bar Y)e_j,e_i\rangle,\qquad
 T=\sum_{i,j=0}^7C_{ij}\partial_{X_i}\partial_{Y_j}.
\]
For fixed $z=X\bar Y$, the fiber through $(X,Y)$ is
$(zY',Y')$; its spherical Laplacian on the restriction of a
bihomogeneous polynomial of total degree $m=k+l$ is
\[
 L_{\VV}P=m(m+6)P-2TP.
\]
The identity follows by expanding the Euclidean Laplacian in the fiber
variable $Y'$ while keeping $z$ fixed, then restricting to $|Y'|=1$.
For octonionic right multiplication $R_a(q)=qa$ we have
$R_a^*=R_{\bar a}$ and
\[
 R_a^*R_b+R_b^*R_a=2\langle a,b\rangle I.
\]
Consequentlx the symmetric part of the linear map
$X\mapsto(X\bar Y)e_j$ is $Y_jI$.
Since $C_{ij}$ is linear in $X$, differentiation gives
\[
 [\Delta_X,T]=2E_Y\Delta_X.
\]
The polarized alternative-algebra identity
\[
 (X\bar e_s)e_j+(X\bar e_j)e_s=2\delta_{js}X
\]
likewise gives $[\Delta_Y,T]=2E_X\Delta_Y$.
Thus if $P$ is separately Euclidean hharmonic, so is $TP$.
The operator $T$ preserves bihomogeneous degrees, since its coefficients
have bidegree $(1,1)$ and it differentiates once in each variable.
Hence $TP\in\PP_{k,l}$, proving the assertion.
\end{proof}

\begin{lemma}\label{bif:spec}
The operator $L_{\VV}$ preserves $\PP_{k,l}$. Its eigenvalues there are among $j(j+6)$, with $j=k+l,k+l-2,\ldots,|k-l|$. For the low bidegrees needed below, the multiplicities are
\[
\begin{array}{c|c|c}
(k,l)&j&\text{multiplicity}\\\hline
(1,1)&0,2&8,56\\
(2,1)\text{ or }(1,2)&1,3&56,224\\
(2,2)&0,2,4&35,350,840.
\end{array}
\]
\end{lemma}
%welches ergebnis stimmt, rechnjngen beide kontrollieren
\begin{proof} The fiber through $(X,Y)$ is $Y'\mapsto(zY',Y')$ with $z=X\bar Y$ and is a round unit $\Sph^7$. Restriction of a product of harmonics of degrees $(k,l)$ to the fiber decomposes into degrees $k+l,k+l-2,\ldots,|k-l|$. The addition theorem gives the vertical trace
\begin{equation}\label{bif:trace}
\operatorname{tr}(L_{\VV}|_{\PP_{k,l}})=\dim\PP_{k,l}\,[k(k+6)+l(l+6)].
\end{equation}
To justify the trace formula, let $d_k=\dim\HH_k(\Sph^7)$ and choose orthonormal bases of homogeneous harmonics $\phi_\alpha(X)$ and $\psi_\beta(Y)$. The addition theorem and homogeneity imply, on the unit spheres,
\[
\sum_\alpha\phi_\alpha(X)\partial_{X_i}\phi_\alpha(X)=kd_kX_i,
\qquad
\sum_\beta\psi_\beta(Y)\partial_{Y_j}\psi_\beta(Y)=ld_lY_j.
\]
For the mixed operator $T$ introduced in Lemma~\ref{bif:bigradeproof}, its diagonal trace therefore equals
\[
\operatorname{tr}T=kld_kd_l\int\langle(X\bar Y)Y,X\rangle=kld_kd_l.
\]
Thus $\operatorname{tr}L_{\VV}=[(k+l)(k+l+6)-2kl]d_kd_l$, as claimed in \eqref{bif:trace}. 
For $(2,2)$ the fiber-constant part is exactly $\{H(X\bar Y):H\in\HH_2(\Sph^7)\}$, of dimension $35$: for fixed $Y$ the map $X\mapsto X\bar Y$ is orthogonal, so any fiber-constant $(2,2)$ function is an eigenfunction of $\Delta_{\Sph^7_z}$ with eigenvalue $16$. The multiplicity calculation uses three separate ingredients: preservation of each bidegree, the list of possible fibrewise degrees, and the trace identity \eqref{bif:trace}. Preservation was established in Lemma~\ref{bif:bigradeproof}. For $(2,2)$ the three eigenvalue multiplicities follow already from the dimension, trace and fiber-constant multiplicity: if $d_2+d_4=1190$, then $16d_2+40d_4=1225\cdot32$, whence $d_2=350$, $d_4=840$. For $(2,1)$, the analogous equations $d_1+d_3=280$ and $7d_1+27d_3=280\cdot23$ yield $d_1=56$, $d_3=224$. The restrictions to the fibers are products of two spherical harmonics after an orthogonal change of variable; the standard harmonic-product decomposition gives only degrees $k+l,k+l-2,\ldots,|k-l|$.
\end{proof}

\begin{lemma}\label{bif:high}
At $t=27$, every bidegree with $k+l\ge5$ has strictly positive
Jacobi eigenvalues.
\end{lemma}
\begin{proof}
Put $m=k+l$. Since $j\le m$ and $k^2+l^2=\frac12(m^2+(k-l)^2)$, the eigenvalue formula
\eqref{bif:eigen} gives
\[
\mu_{k,l,j}\ge(k-l)^2+6m+\frac{m(m+6)}{27}-28.
\]
For $m\ge5$ the right-hand side is strictly positive, as is seen at
$m=5$ and by monotonicity in $m$.
\end{proof}

At $t=27$ the eigenvalue on a vertical-degree-$j$ component is
\begin{equation}\label{bif:eigen}
\mu_{k,l,j}=2[k(k+6)+l(l+6)]-\frac{26}{27}j(j+6)-28.
\end{equation}
Lemma~\ref{bif:high} excludes all bidegrees with $k+l\ge5$.
Thus only $m\le4$ matters. In the following count, each listed multiplicity is the dimension of the corresponding vertical eigenspace, summed over both ordered bidegrees when two are displayed. The complete nonpositive count is
\begin{center}\begin{tabular}{@{}lllrr@{}}\toprule
$(k,l)$&$j$&sign&negative dimension&zero dimension\\\midrule
$(0,0)$&0&$-$&1&0\\
$(1,0),(0,1)$&1&$-$&16&0\\
$(2,0),(0,2)$&2&$-$&70&0\\
$(1,1)$&0&$0$&0&8\\
$(1,1)$&2&$-$&56&0\\
$(3,0),(0,3)$&3&$0$&0&224\\
$(2,1),(1,2)$&3&$-$&448&0\\
$(2,2)$&4&$-$&840&0\\\midrule
Total&&&1431&232\\\bottomrule
\end{tabular}\end{center}
In particular, using Lemma \ref{bif:bigradeproof} and the low-degree multiplicities of Lemma \ref{bif:spec},
\begin{equation}\label{bif:base}
\operatorname{ind}(\Sigma_{27})=1431,\qquad \operatorname{nul}(\Sigma_{27})=232.
\end{equation}
\subsection{Existence of a nontrivial branch}
The stabilizer $G\cong\operatorname{Spin}(7)$ of an axis $e$ in the first half-spin representation acts isometrically on $g_t$ and preserves $\Sigma$. Set $x=\av{X,e}$ and
\begin{equation}\label{bif:pdef}
p(x)=x^3-\frac3{10}x,\qquad \Delta_Xp=-27p.
\end{equation}
The $G$-fixed subspace of the $232$-dimensional kernel is $\mathbb Rp$: the $224$-dimensional critical subspace splits into the first-half-spin cubic harmonics, with one zonal invariant, and the second-half-spin cubic harmonics, with no invariant; the eight permanent Killing modes have no $G$-invariant vector.

\begin{lemma}\label{bif:crdetails}
Fix $0<\alpha<1$ and set $X_G=C^{2,\alpha}(\Sigma)^G$ and
$Y_G=C^{0,\alpha}(\Sigma)^G$. There is a smooth local curve of
nonzero $G$-invariant normal graphs of minimal hypersurfaces
$(t(\varepsilon),u_\varepsilon)$, with
$t(0)=27$, $u_\varepsilon=\varepsilon p+o(\varepsilon)$ in $X_G$.
All nontrivial solutions in a sufficiently small neighborhood in
$\mathbb R\times X_G$ belong to this curve, up to reparametrisation.
\end{lemma}
\begin{proof}
Use the first variation of the graph area to define the Euler--Lagrange
map $\mathcal F(t,u)$ with respect to the \emph{fixed} smooth reference
measure on $\Sigma$. Multiplication by the inverse of the positive
scalar area prefactor $C_t$ does not alter its zeros. The exact area
integrand is smooth for $t>0$ and small $C^{1}$ graphs; its principal
symbol is positive definite near $(27,0)$. Consequently
$\mathcal F:\mathbb R\times X_G\to Y_G$ is smooth and its
linearization $L_t=D_u\mathcal F(t,0)$ is a self-adjoint uniformly
elliptic operator, with the normalization $L_t=J_t$.
The full operator $J_{27}:C^{2,\alpha}\to C^{0,\alpha}$ is Fredholm
of index zero. Since it commutes with the compact group $G$, its
restriction to the invariant spaces is Fredholm of index zero: one
may restrict a $G$-equivariant elliptic parametrix, obtained by
averaging any parametrix over $G$. The spectral decomposition and
Lemma~\ref{bif:repr} give $\ker(L_{27}|_{X_G})=\mathbb Rp$.
Self-adjointness and elliptic regularity then identify the range as
\[
\operatorname{ran}(L_{27}|_{X_G})
=\{f\in Y_G:\int_\Sigma fp\,d\mu=0\}.
\]
Finally, the eigenvalue of $p$ on the trivial branch is
$\mu(t)=27/t-1$, so
\[
\int_\Sigma p\,D_{tu}\mathcal F(27,0)p\,d\mu
=-\frac1{27}\int_\Sigma p^2\,d\mu\ne0.
\]
This is the transversality condition. The simple-eigenvalue
bifurcation theorem yields the asserted curve and local description.
For each sufficiently small solution, the normal-graph equation is a
smooth uniformly elliptic quasilinear equation. Elliptic Schauder
regularity and bootstrapping therefore imply that the graph function
is smooth on $\Sigma$.
\end{proof}

\paragraph{The branch expansion.}
The second-order correction is
\begin{equation}\label{bif:wdef}
w_2=-\frac{27}{10}x^6+\frac{81}{455}x^4
+\frac{39609}{40040}x^2-\frac{15957}{200200}.
\end{equation}
It satisfies, as an exact polynomial identity,
\begin{equation}\label{bif:weq}
\left(\frac{28}{27}(-\Delta_X)-28\right)w_2
=-2(54p^2-|\nabla_Xp|^2).
\end{equation}
The normalized moments are
\begin{equation}\label{bif:moments}
\av{p^2}=\frac7{1600},\quad
\av{|\nabla p|^2}=\frac{189}{1600},\quad
\av{w_2}=\frac{27}{3200}.
\end{equation}
To identify the reduced quartic coefficient, write the normalized
area expansion at $t=27$ as $E_2(u)+E_3(u)+E_4(u)+O(\|u\|^5)$,
with the homogeneous terms specified in \eqref{bif:quartic}.
The order-$\varepsilon^2$ equation is
$J_0w_2+DE_3(p)=0$ on the complement of the critical mode.
Consequently, substitution of
$u=\varepsilon p+\varepsilon^2w_2+O(\varepsilon^3)$ gives
\[
K=E_4(p)+DE_3(p)[w_2]+\tfrac12\langle w_2,J_0w_2\rangle
 =E_4(p)-\tfrac12\langle w_2,J_0w_2\rangle.
\]
Here $E_4(p)$ is the direct quartic contribution, whereas the second
term accounts for the eliminated noncritical modes. Both are evaluated
using the exact spherical moments recorded in Section~\ref{bif:momentappendix}.
The reduced quartic coefficient and branch parameter are
\begin{equation}\label{bif:kc}
K=-\frac{4248690379}{230630400000},\qquad
c=\frac{12746071137}{28028000},\qquad
K=-\frac{c}{108}\av{p^2}.
\end{equation}
In particular the amplitude direction contributes one negative eigenvalue, since its effective second variation is $8K\varepsilon^2+o(\varepsilon^2)<0$. The parameter curve is even in the amplitude. Indeed, choose an element of $\operatorname{SO}(8)$ taking $e$ to $-e$. Its induced action on normal graphs sends the zonal cubic $p(x)$ to $p(-x)=-p(x)$ and preserves the graph equation and the parameter $t$. Parametrize the local bifurcating curve by its coefficient along $p$, using the $L^2$-projection onto $\mathbb Rp$. The transformed solution has amplitude $-\varepsilon$. Local uniqueness of the Crandall--Rabinowitz curve with this amplitude normalization then gives $t(-\varepsilon)=t(\varepsilon)$. Since the curve is smooth, its Taylor expansion has no odd powers, and hence the parameter remainder is $O(\varepsilon^4)$.
\subsection{The effective Hessian}
Set $a=28/27$ and $\beta=2(1/27-1)$. Writing
$A=|\nabla_Xu|^2$, $B=|\nabla_Yu|^2$, $C=\av{\bar z\nabla_Xu,\nabla_Yu}$, and $D=a(A+B)+\beta C$, the quartic Taylor polynomial of $\mathcal A_t/C_t$ at $(t,u)=(27,0)$ is
\begin{align}\label{bif:quartic}
\int\bigg[1-14u^2+\frac12D+2u(A-B)
+\frac{266}{3}u^4+(4-7a)u^2(A+B)\\
\hspace{22mm}-7\beta u^2C-\frac18D^2\bigg]+O(5).\nonumber
\end{align}
For the critical space $E$ let $P$ denote the $L^2$-projection and $J_0=J_{27}$. The restriction of $J_0$ to $E^\perp$ has no zero eigenvalue,
by the spectral decomposition above; its inverse is understood on
$E^\perp$, for instance as the elliptic inverse between the corresponding
Sobolev spaces. Eliminating the $E^\perp$ component of a variation
produces the Schur complement on the critical space.
The effective order-$\varepsilon^2$ operator is
\begin{equation}\label{bif:schur}
B_{\mathrm{eff}}=PJ_2P-PJ_1(J_0|_{E^\perp})^{-1}J_1P.
\end{equation}
For $v=v(X)\in\HH_3(\Sph^7_X)$,
\begin{equation}\label{bif:jone}
J_1v=4(54pv-\av{\nabla p,\nabla v}).
\end{equation}
For $v=v(Y)\in\HH_3(\Sph^7_Y)$, $J_1v=0$; a proof of this operator identity is given below.

For $v=v(X)$ the direct coefficient of $J_2$, divided by $\av{v^2}$, is
\begin{align}\label{bif:direct}
Q(v)=\frac1{\av{v^2}}\av{&\left(4w_2+(8-14a)p^2-\frac{a^2}{2}|\nabla p|^2+\frac{c}{729}\right)|\nabla v|^2
-a^2\av{\nabla p,\nabla v}^2\\
&+8v\av{\nabla w_2,\nabla v}+(32-56a)pv\av{\nabla p,\nabla v}
+\left(1064p^2+(8-14a)|\nabla p|^2\right)v^2}.\nonumber
\end{align}
Representative functions for the $\mathbf{27}$, $\mathbf{77}$ and Killing $\mathbf7$ components are
\[
v_{27}=x_1(x_2^2-x_3^2),\qquad v_{77}=x_2x_3x_4,\qquad
v_7=(3x_1^2-3/10)x_2,
\]
where $x_1,x_2,x_3,x_4$ are mutually orthogonal coordinate functions on the first $\Sph^7$.
The exact computations give
\begin{center}\begin{tabular}{@{}lrrr@{}}\toprule
block&$Q$&Schur subtraction&effective coefficient\\\midrule
$\mathbf7$&$16471593/560560$&$16471593/560560$&$0$\\
$\mathbf{27}$&$584875393/20270250$&$16957431/2802800$&$25885262453/1135134000$\\
$\mathbf{77}$&$45822643/3234000$&$1023111/254800$&$106720261/10510500$\\\bottomrule
\end{tabular}\end{center}
In each row the effective coefficient is the direct quadratic coefficient $Q$ minus the Schur subtraction. The subtraction is the quadratic form of $PJ_1(J_0|_{E^\perp})^{-1}J_1P$ on the indicated representative, divided by its squared $L^2$ norm. The first row also follows from the Killing-field identity.

For the second-half-spin cubic block, the following is the \emph{direct} second-order Hessian coefficient, before Schur complementation:
\begin{align}\label{bif:b112}
b_{112}={}&27\left[(8-14a)\av{p^2}-4\av{w_2}-\frac{a^2}{2}\av{|\nabla p|^2}+\frac c{729}\right]\\
&+1064\av{p^2}+(8-14a)\av{|\nabla p|^2}
-\left(\frac{26}{27}\right)^2\frac{729}{1600}
=\frac{2132701037}{126126000}>0.\nonumber
\end{align}
The potentially missing mixed term $-7\beta u^2C$ does have an $\varepsilon^2s^2$ coefficient. It integrates to zero because, for fixed $X$, the vector field $Y\mapsto(Y\bar X)\nabla_Xp$ is Killing on $\Sph^7_Y$ and $\int vK_Xv=0$.

\begin{lemma}\label{bif:112schur}
For $p=p(X)\in\HH_3(\Sph^7_X)$ and $v=v(Y)\in\HH_3(\Sph^7_Y)$,
$J_1v=0$ as an operator identity. Consequently the Schur correction on
$\HH_3(V_-)$ vanishes, and the effective coefficient equals the direct
coefficient in \eqref{bif:b112}.
\end{lemma}
\begin{proof}
The cubic part of the normalized area functional is
\[
 E_3(u)=2\int u\bigl(|\nabla_Xu|^2-|\nabla_Yu|^2\bigr).
\]
For an arbitrary smooth test function $w=w(X,Y)$, since $\nabla_Yp=0$
and $\nabla_Xv=0$, its third differential satisfies
\[
 D^3E_3(0)[p,v,w]
 =4\int\bigl(v\langle\nabla_Xp,\nabla_Xw\rangle
                 -p\langle\nabla_Yv,\nabla_Yw\rangle\bigr).
\]
Integrating by parts in the two sphere variables, and using
$\Delta_Xp=-27p$ and $\Delta_Yv=-27v$, gives
\[
 D^3E_3(0)[p,v,w]
 =-4\int v(\Delta_Xp)w+4\int p(\Delta_Yv)w
 =108\int pvw-108\int pvw=0.
\]
As $w$ is arbitrary, $J_1v=0$. In particular
$PJ_1(J_0|_{E^\perp})^{-1}J_1Pv=0$, and the claimed equality follows.
\end{proof}

\section{Representation-theoretic reduction and index count}
The critical eigenspace decomposes into irreducible $\operatorname{Spin}(7)$-modules. The effective Hessian and the Jacobi fields induced by ambient isometries determine the index and nullity along the bifurcating branch.
\begin{lemma}\label{bif:repr}
Let $G=\operatorname{Spin}(7)$ be the stabilizer of a unit vector in one of the two real half-spin modules $V_+$ of $\operatorname{Spin}(8)$. Under the standard triality embeddings, the other half-spin module $V_-$ and the vector module $V_0$ restrict to the eight-dimensional spin module of $G$, while $V_+|_G=\mathbb R\oplus\mathbb R^7$. Consequently
\[
\HH_3(V_+)|_G=\mathbf1\oplus\mathbf7\oplus\mathbf{27}\oplus\mathbf{77},
\qquad \HH_3(V_-)|_G=\mathbf{112}.
\]
Each displayed summand is irreducible and no two are isomorphic. The eight-dimensional fiber-constant kernel is a copy of the spin module and has no $G$-fixed vector.
\end{lemma}
\begin{proof}
The stabilizer of a unit vector in the eight-dimensional spin representation of $\operatorname{Spin}(8)$ is $\operatorname{Spin}(7)$. The restrictions of the other two eight-dimensional triality modules are the irreducible eight-dimensional spin module. Write $V_+=\mathbb Re\oplus W$, where $W$ is the standard seven-dimensional $G$-module. A cubic harmonic polynomial on $\mathbb R e\oplus W$ is determined uniquely by its harmonic coefficients of degrees $0,1,2,3$ in $W$: solve the ambient harmonicity equation successively for the powers of the coordinate along $e$. This yields a $G$-equivariant isomorphism
\[
\HH_3(\mathbb R\oplus W)\simeq
\HH_0(W)\oplus\HH_1(W)\oplus\HH_2(W)\oplus\HH_3(W).
\]
Their dimensions are $1,7,27,77$. Each is irreducible under $\operatorname{SO}(7)$, hence under its spin cover. For the second assertion, complexify the eight-dimensional spin module $S=V_-\otimes\mathbb C$ of the compact group of type $B_3$. Its highest weight is $\omega_3=(\tfrac12,\tfrac12,\tfrac12)$ in orthogonal coordinates. The cube of a highest-weight vector is a highest-weight vector of weight $3\omega_3$, so the irreducible module $V_{3\omega_3}$ occurs in $\operatorname{Sym}^3 S$. The Weyl dimension formula for $B_3$, with $\rho=(\tfrac52,\tfrac32,\tfrac12)$, gives
\[
 \dim V_{3\omega_3}
 =\frac{4\cdot3\cdot2}{(5/2)(3/2)(1/2)}
 \frac{16-9}{25/4-9/4}
 \frac{16-4}{25/4-1/4}
 \frac{9-4}{9/4-1/4}=112.
\]
The equivariant trace map $\operatorname{Sym}^3 S\to S$ is onto (the metric contraction is nonzero), so its kernel, the complexified harmonic cubics, has dimension $\binom{10}{3}-8=112$. Since $V_{3\omega_3}$ cannot map to $S$ and already has this full dimension, the harmonic cubic kernel equals $V_{3\omega_3}$ and is irreducible. This proves the asserted real irreducibility as well, since any proper real invariant subspace would complexify to a proper complex invariant subspace. Finally, the degree-one harmonics of the quotient variable $X\bar Y$ form an eight-dimensional spin module under $G$, so they have no fixed vector.
\end{proof}
\begin{lemma}\label{bif:orbit}
For every sufficiently small nonzero member of the bifurcating branch, the
identity component of its stabilizer in $K=\operatorname{Spin}(9)$ is
$G=\operatorname{Spin}(7)$. Consequently the normal components of ambient
Killing fields span a $15$-dimensional subspace of the Jacobi kernel.
\end{lemma}
\begin{proof}
Choose an $H$-invariant tubular neighborhood of $\Sigma$ and identify
nearby unparametrized $C^{2,\alpha}$ hypersurfaces with their normal graph
functions in $\mathcal X=C^{2,\alpha}(\Sigma)$. The identity component of
the stabilizer of $\Sigma$ in $K$ is $H=\operatorname{Spin}(8)$; the
infinitesimal normal action of $\mathfrak k/\mathfrak h$ is an injective
map onto the eight-dimensional Killing space $W\subset\mathcal X$.
To see injectivity, a Killing field with zero normal component along
$\Sigma$ is tangent to $\Sigma$ and its flow preserves $\Sigma$;
hence it belongs to the Lie algebra of its stabilizer $H$. Thus the
normal action of $\mathfrak k/\mathfrak h$ has zero kernel. Its
dimension is $36-28=8$. The normal action can be identified explicitly, without relying on
an abstract equivalence of triality modules. Write the octonionic Hopf
map in unit-sphere coordinates as
\[
 \pi(a,b)=(s,z)=(|a|^2-|b|^2,2a\bar b)\in S^8\subset\mathbb R\oplus\mathbb O.
\]
The standard $K=\operatorname{Spin}(9)$ action on the Hopf fibration
covers the rotation action on $S^8$. For each $A\in\mathbb O$, take the
infinitesimal base rotation with $\dot s=\langle z,A\rangle$ and
$\dot z=-sA$. Its lift is an ambient Killing field. On the balanced
Clifford hypersurface, $s=0$ and $z=X\bar Y$. Along the normal unit
variation $(\cos\theta X,\sin\theta Y)$ one has
$s=\cos(2\theta)$ and $\partial_\theta s|_{\theta=\pi/4}=-2$.
Since the lifted field is projectable, its normal component is therefore
\[
 k_A(X,Y)=-\tfrac12\langle X\bar Y,A\rangle.
\]
These eight linearly independent functions span precisely the degree-one
fiber-constant space $H_1(X\bar Y)$. The eight rotations above represent
$\mathfrak k/\mathfrak h$, hence $W$ is the full infinitesimal normal
orbit space, not merely a chosen subspace of Killing variations.
Let $P_W$ be the $L^2$-orthogonal projection onto $W$ and set
$\mathcal S=\ker P_W$. The $H$-action on $\mathcal X$ is linear and
preserves $\mathcal S$.

For the local uniqueness assertion, choose
an $\operatorname{Ad}(H)$-invariant complement $\mathfrak m$ to
$\mathfrak h$ in $\mathfrak k$. For $A\in\mathfrak m$ and
$u\in\mathcal S$ near zero, let $\Phi(A,u)$ be the normal graph of
$\exp(A)\cdot\operatorname{graph}(u)$ and consider
\[
 F(A,u)=P_W\Phi(A,u).
\]
At $(0,0)$, $D_AF:\mathfrak m\to W$ is an isomorphism: it is precisely
the normal infinitesimal orbit map. The Banach-space implicit function
theorem therefore gives neighborhoods in which
\[
 F(A,u)=0\quad\Longrightarrow\quad A=0,
 \qquad u\in\mathcal S.
\]
%begründung einfügen
(The solution $A=0$ exists because $u\in\mathcal S$.) Consequently,
if $u\in\mathcal S$ is sufficiently small and $k\in K$ stabilizes its
graph and is sufficiently close to the identity, write $k=\exp(A)h$
with $A\in\mathfrak m$ small and $h\in H$ near the identity. Since
$h\cdot u\in\mathcal S$ and the final graph is again $u\in\mathcal S$,
uniqueness forces $A=0$, hence $k\in H$. This argument needs only the
local action and does not assume a global slice theorem for embeddings.

The graph $u_\varepsilon$ is $G$-fixed. Since $W$ restricts to the
irreducible eight-dimensional spin representation of $G$, it has no
$G$-fixed vector. As $P_W$ is $G$-equivariant,
$P_Wu_\varepsilon=0$; thus the branch belongs to $\mathcal S$.
The preceding local uniqueness shows that the Lie algebra of its
$K$-stabilizer lies in $\mathfrak h$.

Let $E_+=\mathcal H_3(V_+)$, viewed as an $H$-submodule of the fixed
graph space, and denote its $H$-equivariant orthogonal projection by
$\Pi_+$. Since $u_\varepsilon$ is $G$-invariant and
$E_+^G=\mathbb Rp$, we have
\[
 \Pi_+u_\varepsilon=a_\varepsilon p,
 \qquad a_\varepsilon=\varepsilon+O(\varepsilon^2)\ne0
\]
for sufficiently small $\varepsilon\ne0$. Every element of $H$
stabilizing $u_\varepsilon$ must therefore stabilize $p$.
The identity component of the stabilizer of the nonzero zonal cubic
$p(x)=x^3-3x/10$ in $H$ is the axis stabilizer $G$. Indeed, on the
unit sphere of $V_+$ the function $p(\langle X,e\rangle)$ has a unique
global maximum at $X=e$: its endpoint values are $p(1)=7/10$ and
$p(-1)=-7/10$, while the only interior critical points satisfy
$x^2=1/10$ and have values of absolute size $1/(5\sqrt{10})<7/10$.
Any element preserving $p$ must therefore fix $e$. Conversely, the
axis stabilizer fixes $p$, so its identity component is precisely $G$.
Since $G$ fixes $u_\varepsilon$ by construction, the connected
stabilizer is exactly $G$, and the dimension of the $K$-orbit is
$\dim K-\dim G=36-21=15$. Normal components of its infinitesimal
actions are Jacobi fields because the hypersurfaces are minimal.
\end{proof}

\begin{lemma}\label{bif:perturb}
Let $\mathfrak q_\varepsilon$ be the second variation of the normalized graph area
at $(t(\varepsilon),u_\varepsilon)$, evaluated on functions on the fixed
reference hypersurface and integrated against its fixed reference measure.
Assume the smooth graph branch and the quartic expansion above. Then
$\mathfrak q_\varepsilon$ is a $C^2$ family of bounded symmetric forms on
$H^1(\Sigma)$, with expansion
\[
\mathfrak q_\varepsilon=\mathfrak q_0+\varepsilon\mathfrak q_1+
\varepsilon^2\mathfrak q_2+o(\varepsilon^2)
\quad\text{in }\mathcal L(H^1,H^{-1}).
\]
Let $E=\ker J_0$, $F=E^{\perp_{L^2}}$, and let $P,Q$ be the corresponding
projections. If $P J_1 P=0$, the reduced form on $E$ has leading
coefficient
\[
B_{\rm eff}=PJ_2P-PJ_1Q(QJ_0Q|_F)^{-1}QJ_1P.
\]
Here $J_i$ denote the distributional operators representing the
coefficients $\mathfrak q_i$ with respect to the \emph{same fixed}
reference measure. The complement contributes its unchanged index and
has no zero modes for sufficiently small $\varepsilon$.
\end{lemma}
\begin{proof}
The exact graph integrand is smooth in $t$, $u$, and the first
derivatives of $u$ in a neighborhood of the reference graph. The
bifurcating curve depends smoothly on $\varepsilon$ in
$C^{2,\alpha}$, and $t(\varepsilon)$ is smooth. Hence the
coefficients of its second variation, written in divergence form
relative to the fixed reference measure, depend twice continuously
differentiably on $\varepsilon$ in $C^0$. Uniform ellipticity persists
for small $\varepsilon$, and the coefficient expansion gives the
asserted expansion in $\mathcal L(H^1,H^{-1})$. The self-adjoint elliptic
operator $J_0$ has compact resolvent. Its restriction to $F$ has a
spectral gap at zero. Write $H^1_F=H^1(\Sigma)\cap F$ and let
$(H^1_F)^*$ denote its continuous dual, with the $L^2$ pairing used
to identify distributions with functionals. The restriction of the
quadratic form of $J_0$ therefore induces a bounded isomorphism
$J_0:H^1_F\to(H^1_F)^*$. In particular, the inverse in the Schur
complement acts on the restricted functional $z\mapsto\langle J_1v,z\rangle$
for $z\in H^1_F$.
For each $v\in E$ the complementary equation
$\mathfrak q_\varepsilon(v+w,z)=0$ for every $z\in H^1\cap F$
has a unique solution $w=w_\varepsilon(v)\in H^1\cap F$ near zero,
by the bounded inverse theorem, and
\[
w_\varepsilon(v)=-\varepsilon(J_0|_F)^{-1}QJ_1v
+O(\varepsilon^2)\|v\|.
\]
Substitution into the form gives the asserted Schur coefficient and a
remainder $o(\varepsilon^2)\|v\|^2$, uniformly on the unit sphere of
the finite-dimensional space $E$. Completing the square by the
invertible complementary block is a congruence of symmetric forms,
so preserves index and nullity. The complementary block retains its
inertia by the spectral gap and continuity of the elliptic forms.
Finally, the fixed graph Hessian has the same index and nullity as the
geometric Jacobi form: they are the same second variation expressed
using different positive densities and normal variation coordinates.
\end{proof}

The hypothesis $PJ_1P=0$ in Lemma~\ref{bif:perturb} holds. Indeed the cubic graph energy gives $J_1v=4(54pv-\langle\nabla_Xp,\nabla_Xv\rangle)$ for $v$ in the first cubic block, which is even in $X$, hence orthogonal to both cubic blocks and to the eight-dimensional mixed Killing block by bidegree parity. On the second cubic block $J_1v=0$ by Lemma~\ref{bif:112schur}. On the mixed Killing block, choose smooth normal Killing variations $k_\varepsilon$ with $J_\varepsilon k_\varepsilon=0$ for the fixed-form operators. Differentiating gives $J_1k_0+J_0\dot k_0=0$; applying $P$ gives $PJ_1k_0=0$. The eight values $k_0$ span the mixed Killing block. Therefore $PJ_1P=0$.

By Lemma~\ref{bif:repr}, equivariance and self-adjointness make the effective Hessian scalar on each irreducible summand. The coefficient on the trivial summand $\mathbf1$ is negative by \eqref{bif:kc}; those on $\mathbf{27}$, $\mathbf{77}$, and $\mathbf{112}$ are positive by the calculations in Sections~5 and~7.

The exact Killing--Jacobi fields account for both zero blocks. The eight normal Killing fields already present on the reference hypersurface span $\mathbf8$. For the other seven directions, the infinitesimal action of $\mathfrak{spin}(8)/\mathfrak{spin}(7)$ on the zonal cubic $p$ is injective: its kernel is precisely the Lie algebra of the stabilizer of $p$, by Lemma~\ref{bif:orbit}. The image is a seven-dimensional $G$-submodule of $\mathcal H_3(V_+)$ and hence equals its unique $\mathbf7$ summand by Lemma~\ref{bif:repr}. Consequently these seven orbit directions have independent leading terms in $\mathbf7$. Together with the eight permanent modes, they give fifteen linearly independent exact Jacobi fields for every sufficiently small nonzero branch parameter.

The critical space and the signs of the reduced form are
\[
\begin{array}{c|rrrrrr}
\text{summand}&\mathbf1&\mathbf7&\mathbf{27}&\mathbf{77}&\mathbf{112}&\mathbf8\\\hline
\text{dimension}&1&7&27&77&112&8\\
\text{reduced sign}&-&0&+&+&+&0
\end{array}
\]
The zero entries alone would not imply exact nullity; that conclusion uses the Killing fields just described. The spectral gap on $E^\perp$ and Lemma~\ref{bif:perturb} exclude additional zero modes, while the negative scalar direction contributes one to the index. Thus, for all sufficiently small nonzero $\varepsilon$,
\[
\operatorname{ind}=1431+1=1432,\qquad \operatorname{nul}=8+7=15.
\]

\begin{remark}
Bifurcation of Clifford hypersurfaces and tori in other geometries is
known; in particular, related constructions concern constant-mean-curvature
Clifford tori in Berger spheres and Willmore Clifford tori in Berger
three-spheres. For Willmore Clifford tori in Berger three-spheres, see Rodrigues~\cite{Rodrigues2026}; for Jacobi eigenvalue bounds and rigidity of surfaces in contracted Berger three-spheres, see Batista--Mendes~\cite{BatistaMendes2026}. These results concern different dimensions, ambient geometries, or variational problems from the octonionic minimal-hypersurface bifurcation considered here.
\end{remark}
%rechnung zu beginn des abschnitts nochmal kontrolle

\section{Moment calculations for the critical blocks}\label{bif:momentappendix}
The spherical moments and projection norms entering the effective Hessian are evaluated by polynomial identities and exact rational integration. All
averages in this section are normalized spherical averages on $S^7$.
For orthonormal coordinate functions $x_1,\ldots,x_r$ and nonnegative integers
$\alpha_i$, rotational symmetry and the beta integral give
\begin{equation}\label{bif:momentformula}
 \av{\prod_{i=1}^r x_i^{2\alpha_i}}
 =\frac{\prod_{i=1}^r(\frac12)_{\alpha_i}}
 {(4)_{\alpha_1+\cdots+\alpha_r}},\qquad
 \av{\prod_{i=1}^r x_i^{a_i}}=0
 \quad\text{if some $a_i$ is odd}.
\end{equation}
Here $(a)_j=a(a+1)\cdots(a+j-1)$ is the rising factorial.
For polynomial restrictions to $S^7$, the tangential gradient pairing and
positive spherical Laplacian are computed algebraically by
\begin{equation}\label{bif:momentoperators}
 \langle\nabla_Sf,\nabla_Sg\rangle
 =\sum_i\partial_if\partial_ig-(Ef)(Eg),\qquad
 Lf=-\sum_i\partial_i^2f+E(E+6)f,
\end{equation}
where $E=\sum_i x_i\partial_i$ acts on the ambient polynomial. In the
second expression $E(E+6)f$ means $E(Ef+6f)$.
The expressions in \eqref{bif:momentoperators} are valid on restrictions
of arbitrary polynomials, not only homogeneous ones.

For $v$ in the first cubic block set
\[
 f_v=4\bigl(54pv-\langle\nabla_Sp,\nabla_Sv\rangle\bigr).
\]
Let $\pi_jf_v$ be its degree-$j$ spherical harmonic projection.
Since $L\pi_jf_v=j(j+6)\pi_jf_v$, the quantities
$M_i=\av{f_v L^if_v}$ satisfy
\begin{equation}\label{bif:momentvandermonde}
 M_i=\sum_{j\in\mathcal J_v}[j(j+6)]^i\,\av{(\pi_jf_v)^2}
 \qquad(0\le i<|\mathcal J_v|).
\end{equation}
The sets $\mathcal J_v$ are $\{2,4,6\}$ for $\mathbf{27}$ and
$\mathbf7$, and $\{4,6\}$ for $\mathbf{77}$. Indeed $f_v$ is an even
polynomial of degree at most six. Under the stabilizer of the distinguished
axis, degree zero carries only the trivial representation, whereas
$\mathbf{27}$, $\mathbf{77}$ and $\mathbf7$ are nontrivial and inequivalent
to it. In degree two the stabilizer representation contains no
$\mathbf{77}$, giving the asserted smaller set for that block.

Choose mutually orthogonal coordinate functions
$x_1,x_2,x_3,x_4$ on the unit $S^7$, and use the following representatives:
\[
 p=x_1^3-\frac3{10}x_1,\qquad
 v_{27}=x_1(x_2^2-x_3^2),\qquad
 v_{77}=x_2x_3x_4,\qquad
 v_7=\left(3x_1^2-\frac3{10}\right)x_2.
\]
For each representative, form $f_v$ as above, expand all ambient
polynomials, and replace each monomial by its average from
\eqref{bif:momentformula}. To recover the projections without solving
the moment system, one may equivalently apply the spectral polynomial
\[
 \pi_j f_v=\prod_{\substack{k\in\mathcal J_v\\k\ne j}}
 \frac{L-k(k+6)}{j(j+6)-k(k+6)}\,f_v.
\]
The resulting polynomials satisfy $L\pi_jf_v=j(j+6)\pi_jf_v$,
$\sum_{j\in\mathcal J_v}\pi_jf_v=f_v$, and
$\av{(\pi_jf_v)(\pi_kf_v)}=0$ for $j\ne k$.
The projection norms can therefore be computed from the moment system or from these spectral projections.

For reproducibility, the eigenvalues used in the projectors are
$\lambda_2=16$, $\lambda_4=40$, and $\lambda_6=72$.
For example, when $\mathcal J_v=\{2,4,6\}$ the projections are
\[
\pi_2f_v=\frac{(L-40)(L-72)}{(16-40)(16-72)}f_v,\quad
\pi_4f_v=\frac{(L-16)(L-72)}{(40-16)(40-72)}f_v,\quad
\pi_6f_v=\frac{(L-16)(L-40)}{(72-16)(72-40)}f_v.
\]
For $\mathbf{77}$, use the two-eigenvalue formula with
$\mathcal J_v=\{4,6\}$. Expanding these polynomials and applying
\eqref{bif:momentformula} gives the table below; the orthogonality
and reconstruction identities provide independent checks.

The following table lists the \emph{unnormalized} squared projection norms
$\av{(\pi_jf_v)^2}$ obtained by inserting the representative cubics into \eqref{bif:momentformula}--\eqref{bif:momentvandermonde}:
\[
\begin{array}{c|c|ccc}
\text{block}&\av{v^2}&j=2&j=4&j=6\\\hline
\mathbf{27}&1/240&1/20&24299/224000&1287/1280\\
\mathbf{77}&1/960&0&59643/1792000&819/10240\\
\mathbf7&27/1600&6561/8000&4831083/896000&8019/1024
\end{array}
\]
The Jacobi eigenvalue on the degree-$j$ component in this one-factor
sector is $\lambda_j=\frac{28}{27}j(j+6)-28$. Thus the Schur subtraction,
with the same normalization as \eqref{bif:direct}, is exactly
\begin{equation}\label{bif:schurexplicit}
 S(v)=\frac1{\av{v^2}}
 \sum_{j\in\mathcal J_v}
 \frac{\av{(\pi_jf_v)^2}}{\frac{28}{27}j(j+6)-28}.
\end{equation}
Every denominator in this sum is nonzero. Direct substitution gives
\[
 S(v_{27})=\frac{16957431}{2802800},\quad
 S(v_{77})=\frac{1023111}{254800},\quad
 S(v_7)=\frac{16471593}{560560}.
\]
Substituting the same representative cubics into
\eqref{bif:direct} and evaluating monomials by
\eqref{bif:momentformula} gives
\[
 Q(v_{27})=\frac{584875393}{20270250},\quad
 Q(v_{77})=\frac{45822643}{3234000},\quad
 Q(v_7)=\frac{16471593}{560560}.
\]
The direct and Schur terms, normalized by $\av{v^2}$ as in
\eqref{bif:direct} and \eqref{bif:schurexplicit}, combine as follows:
\[
\begin{array}{c|ccc}
\text{block}&Q(v)&S(v)&Q(v)-S(v)\\\hline
\mathbf{27}&\frac{584875393}{20270250}&\frac{16957431}{2802800}&\frac{25885262453}{1135134000}\\
\mathbf{77}&\frac{45822643}{3234000}&\frac{1023111}{254800}&\frac{106720261}{10510500}\\
\mathbf7&\frac{16471593}{560560}&\frac{16471593}{560560}&0
\end{array}
\]
In the $\mathbf7$ block, the equality $Q(v_7)=S(v_7)$ is also forced
geometrically by the exact Killing--Jacobi fields along the branch. 

%am ende nochmals komtrollieren

\bigskip

\noindent
\textsc{Anna Siffert}\\
Universit\"at M\"unster, Mathematisches Institut,
Einsteinstr.\ 62, 48149 M\"unster, Germany\\
\texttt{asiffert@uni-muenster.de}

\end{document}